\documentclass{amsart}

\usepackage[a4paper,hmargin=3.5cm,vmargin=4cm]{geometry}
\usepackage{amsfonts,amssymb,amscd,amstext}
\usepackage{graphicx}
\usepackage[dvips]{epsfig}
\usepackage{mathpazo}
\usepackage{enumerate}

\def\r{\mathbb{R}}
\def\n{\mathbb{N}}
\def\c{\mathbb{C}}

\def\s{\mathbb{S}}
\def\d{\mathbb{D}}

\def\z{\mathbb{Z}}

\def\Xcal{\mathcal{X}}

\def\Kcal{\mathcal{K}}

\def\Hcal{\mathcal{H}}

\def\Ocal{\mathcal{O}}

\def\Tcal{\mathcal{T}}
\def\Ccal{\mathcal{C}}
\def\Ecal{\mathcal{E}}
\def\Scal{\mathcal{S}}

\def\Hgot{\mathfrak{H}}
\def\Ggot{\mathfrak{G}}

\newtheorem{theorem}{Theorem}[section]
\newtheorem{proposition}[theorem]{Proposition}
\newtheorem{claim}[theorem]{Claim}
\newtheorem{lemma}[theorem]{Lemma}
\newtheorem{corollary}[theorem]{Corollary}
\newtheorem{remark}[theorem]{Remark}

\theoremstyle{definition}

\newtheorem{definition}[theorem]{Definition}

\numberwithin{equation}{section}
\numberwithin{figure}{section}

\begin{document}

\title[On harmonic quasiconformal immersions of surfaces in $\r^3$]
{On harmonic quasiconformal immersions of surfaces in $\r^3$}

\author[A.~Alarc\'{o}n]{Antonio Alarc\'{o}n}
\address{Departamento de Geometr\'{\i}a y Topolog\'{\i}a \\
Universidad de Granada \\ E-18071 Granada \\ Spain}
\email{alarcon@ugr.es}

\author[F.J.~L\'{o}pez]{Francisco J. L\'{o}pez}
\address{Departamento de Geometr\'{\i}a y Topolog\'{\i}a \\
Universidad de Granada \\ E-18071 Granada \\ Spain}
\email{fjlopez@ugr.es}


\thanks{Research partially
supported by MCYT-FEDER research project MTM2007-61775 and Junta
de Andaluc\'{i}a Grant P09-FQM-5088}

\subjclass[2010]{53C43; 53C42, 	30F15}
\keywords{Harmonic immersions of Riemann surfaces, quasiconformal mappings, Gauss map}

\begin{abstract}
This paper is devoted to the study of the global properties of harmonically immersed Riemann surfaces in $\r^3.$ We focus on the geometry of complete  harmonic immersions with quasiconformal Gauss map, and in particular, of those with finite total curvature.  We pay special attention to the construction of new examples with significant geometry.
\end{abstract}

\maketitle

\thispagestyle{empty}


\section{Introduction}
A smooth map $X:M \to N$ between Riemannian manifolds $M$ and $N,$ $M$ compact, is said to be {\em harmonic} if it is a critical point of the Dirichlet functional $\int_M \|dX\|^2 dV.$ If $M$ is not compact then $X$ is said to be harmonic if $X|_{\Omega}$ does for all compact region $\Omega\subset M.$ Harmonic maps model the extrema of the energy functionals associated to some physical phenomena in viscosity, dynamic of fluids, electromagnetism, cosmology...  When $N=\r^n,$ they simply correspond to the solutions to the Dirichlet equation $\Delta X=0,$ where $\Delta$ is the Laplace operator on $M.$ See, for instance, the surveys \cite{helein,helein-wood} and the references therein for a good setting.

If  $M$ has dimension two, the Dirichlet functional actually depends only on the conformal structure of the surface. The interplay between conformality, minimality and harmonicity is a core topic in geometrical analysis, and more particularly, in surface theory. It is well known that a conformal immersion of a Riemann surface in the Euclidean three-space is minimal if and only if it is harmonic, and in this case its Gauss map is conformal. However,  harmonic minimal immersions in $\r^3$ are not necessarily  conformal as the following parameterization of the half helicoid  shows $$X:\c\to\r^3, \quad X(z)=\Re(e^z,\imath e^z,\imath z).$$
The questions considered in this paper are all related to the global theory of harmonic immersions of Riemann surfaces in $\r^3.$ This study was initiated by Klotz \cite{K1,K2,K3,K} in the sixties, at the same time when Osserman \cite{O} laid the foundations of the global theory of minimal surfaces. She made an effort to distinguish those facts about minimal surfaces which are special to them from among the many facts which apply to harmonic immersions. However, there have been several difficulties for the development of this topic. Among them,  we emphasize the lack of enlightening examples (other than minimal surfaces) that facilitate the building of intuition and geometric insight, and the fact that Hopf's maximum principle does not apply to this family of immersions.

In this paper we have introduced a Weierstrass type representation for harmonic immersions that permit a better understanding of their geometric and analytic properties, specially those related to the Gauss map. 
Although in general the Gauss map of a non-flat harmonic immersion is far from being conformal, it shares some basic topological properties with the one of minimal surfaces. For instance, it is open and  everywhere regular except at  a discrete set of topological branch points.

For a good understanding of the subsequent results,  the following notations are required. 

Let $M$ be an open Riemann surface, and let $X=(X_j)_{j=1,2,3}:M \to \r^3$ be a harmonic immersion. We denote by $\Ggot:M \to \s^2,$ $\Kcal$ and $\sigma$ the orientation preserving Gauss map, the (nowhere positive) Gauss curvature and the second fundamental form of $X,$ respectively.  The holomorphic 2-form $\Hgot:=\sum_{j=1}^3 (\partial_z X_j)^2$ is called the {\em Hopf differential} of $X,$  where $\partial_z$ denotes  the  complex differential on $M.$ $X$ is conformal (and so minimal) if and only  if $\Hgot=0.$

The immersion $X$ is said to have {\em vertical flux} if $\partial_z X_j$ is exact, $j=1,2.$ Equivalently, $X$ has vertical flux if for any arc-length parameterized closed curve  $\gamma (s) \subset M$ the flux vector $\int_\gamma \mu(s) ds,$  where $\mu(s)$ is  the conormal vector of $X$ at $\gamma(s),$ is vertical. Since $X$ is  harmonic, the flux of $X$ on a curve depends only on its homology class.  

A plane $\Pi$ of $\r^3$ is said to be {\em finite} for $X$ if $X(M)$ and $\Pi$ meet transversally except for a finite number of points, and $X^{-1}(\Pi)$ consists of a finite number of curves.  $X$ is said to be of {\em finite total curvature} (FTC  for short) if $\int_M \|\sigma\|^2 dS<+\infty,$ where $dS$ is the Euclidean area element on $M.$  $X$  is said to be {\em algebraic} if  $M$ is biholomorphic to $\overline{M}-\{Q_1,\ldots,Q_k\},$  where $\overline{M}$ is a compact Riemann surface and $\{Q_1,\ldots,Q_k\} \subset \overline{M},$  and the vectorial holomorphic 1-form  $\partial_z X$ extends meromorphically to  $\overline{M}$ with polar singularities at $Q_j$ for all $j.$ 

The subclass of harmonically immersed Riemann surfaces with {\em quasiconformal} Gauss map (harmonic QC immersions for short) is specially interesting and gives rise to a rich theory. A harmonic $X:M\to\r^3$ is QC if and only if the map $X$ is quasiconformal in the classical sense, that is to say, if the distortion function of $X$
\[
\mathfrak{D}_X:=\frac{\|X_u\|^2+\|X_v\|^2}{2\|X_u\times X_v\|}\quad (z=u+\imath v\text{ conformal parameter on } M)
\]
is bounded (see \cite{Ka} and Remark \ref{re:beltrami}).
In this paper we have proved that this family contains all complete harmonic immersions with FTC. Among other special features,  QC harmonic immersions are quasiminimal in the sense of Osserman \cite{O1}, and support both Huber-Osserman and Jorge-Meeks  type results that generalize the ones for minimal surfaces (see \cite{O,JM}). 
\begin{quote}
{\bf Theorem I.}
{\em A complete harmonic immersion $X$ of an open Riemann surface $M$ into $\r^3$ has FTC if and only if either of the following statements hold:
\begin{itemize}
\item $X$ is algebraic and the Gauss map $\Ggot:M \to \s^2$ of $X$ extends continuously to $\overline{M}.$
\item $X$ is algebraic and QC.
\item $X$ is QC, $X$ has two non parallel  finite planes and either $\sup_M \|\sigma\|^2<+\infty$ or $X$ is proper.   
\end{itemize}
Furthermore, if $X$ has FTC and $M=\overline{M}-\{Q_1,\ldots,Q_k\},$ where $\overline{M}$ is the compactification of $M,$  then
\begin{itemize}
\item $\Ggot:\overline{M} \to \s^2$ is a finitely branched covering and the classical Jorge-Meeks formula
$$\int_M \Kcal dS= -4\pi  Deg(\Ggot)=-2 \pi\big(2 Gen(M)-2+ \sum_{j=1}^k (I_{Q_j}+1)\big),$$ holds, where $Gen(M),$ $Deg(\Ggot)$ and $I_{Q_j}+1$ are  the genus of $\overline{M},$ the number of sheets of $\Ggot$ and the pole order of  $\partial_z X$ at $Q_j$ for all $j,$ respectively. 
\item $X$ is proper and $X(M)$ viewed from infinity looks like a
collection of $k$ flat planes  with multiplicities $I_{Q_j},$ $j=1,\ldots,k,$ that pass through the origin.
\end{itemize}}
\end{quote}

\begin{figure}[ht]
    \begin{center}
    \scalebox{0.5}{\includegraphics{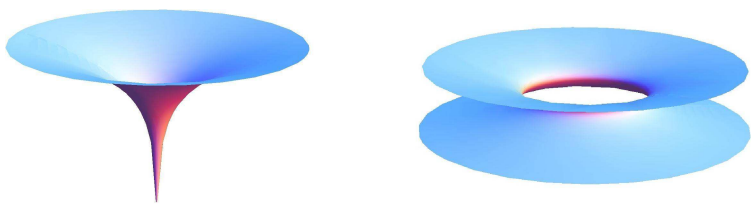}}
        \end{center}
\caption{ The rotational harmonic horn and a rotational harmonic catenoid.}\label{fig:chincheta-catenoide}
\end{figure}

Complete harmonic embeddings of FTC  share some basic topological properties of  embedded minimal surfaces \cite{JM}. 
Outside a compact set in $\r^3,$ a harmonic embedded end of FTC is a graph either with non-zero logarithmic growth ({\em harmonic catenoidal end}) or asymptotic to a plane ({\em harmonic planar end}). Consequently, complete harmonic embeddings of FTC have parallel ends, and ordering them by heights, two consecutive ends have opposite normal. 
In this paper we have classified  complete harmonic embedded annuli with FTC, and in particular,  rotational harmonic immersions. It is worth mentioning that some of the most celebrated classification theorems for embedded minimal surfaces (see \cite{O,JM,lop-ros,schoen}) do not hold in the harmonic framework. 
For instance, there are complete rotational harmonic embedded annuli non-linearly equivalent to a minimal catenoid in both cases with FTC ({\em rotational harmonic catenoids})  and without FTC ({\em the rotational harmonic horn}), see Figure \ref{fig:chincheta-catenoide}. Contrary to the minimal case, one also has:

\begin{quote} 
{\bf Theorem II.}
{\em There exists a complete harmonic embedding $X:M \to \r^3$  of FTC satisfying either of the following properties:
\begin{itemize}
\item  $X(M)$ is a non-rotational annulus in both cases with and without vertical flux (Figure \ref{fig:cat-flux}).
\begin{figure}[ht]
    \begin{center}
    \scalebox{0.3}{\includegraphics{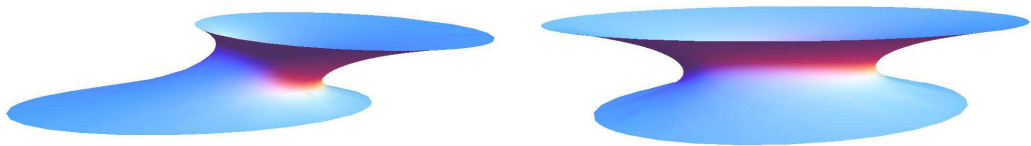}}
        \end{center}
\caption{Two non-rotational harmonic catenoids. On the left: with non-vertical flux. On the right: with vertical flux.}\label{fig:cat-flux}
\end{figure}
\item $X(M)$ is a genus zero surface with vertical flux and more than two ends (Figure \ref{fig:flujo+toro}).
\item $X(M)$ is a twice punctured torus (Figure \ref{fig:flujo+toro}). The existence of this example was pointed out by Weber \cite{we}.
\begin{figure}[ht]
    \begin{center}
    \scalebox{0.5}{\includegraphics{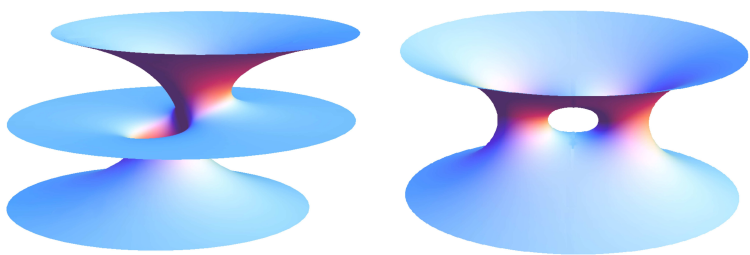}}
        \end{center}
\caption{On the left: a complete harmonic embedding of a three punctured sphere with FTC and vertical flux. On the right:  Weber's surface.}\label{fig:flujo+toro}
\end{figure}
\end{itemize}}
\end{quote}

The study of the Gaussian image, which is one of the fundamental topics in the global theory of minimal surfaces,  can be  naturally extended to  the larger family of harmonic immersions. In this setting, the image of the Gauss map of a complete non-flat harmonic immersion with nowhere vanishing Hopf differential can not omit a neighborhood of an equator, and if in addition the immersion is QC, then it can not lie in $\s^2\cap\{|x_3|<1-\epsilon\},$ $\epsilon>0$ (see Klotz \cite{K}). 
One of the paradigmatic examples in this setting is the rotational harmonic horn, whose Gaussian image coincides, up to reversing the orientation, with the upper open hemisphere minus the North pole. We have improved Klotz's theorems  by dropping the hypothesis about the Hopf differential:
\begin{quote}
{\bf Theorem III.}
{\em Let $X:M\to\r^3$ be a harmonic immersion, where $M$ is an open Riemann surface, and call $\Ggot$ its Gauss map. The following statements hold:
\begin{itemize}
\item If $X$ is complete and  $\Ggot(M) \subset \s^2\cap\{|x_3|>\epsilon\}$ for some $\epsilon>0,$ then $X(M)$ is a plane.
\item If $X$ is complete and QC, and $\Ggot(M) \subset \s^2\cap\{|x_3|<1-\epsilon\}$ for some $\epsilon>0,$ then $X(M)$ is a plane.
\end{itemize}}
\end{quote}

Finally, we have obtained some Privalov's type results for harmonic immersions. A harmonic immersion $X:M \to \r^3$ is said to be {\em asymptotically conformal} if the Beltrami differential of its Gauss map vanishes at the ideal boundary of  $M.$ The following theorem gives some information of how the asymptotic behavior of the Gauss map influences the geometry of the immersion:

\begin{quote}
{\bf Theorem IV.} {\em  Let $X:M\to\r^3$ be a non-flat harmonic immersion, where $M$ is an open Riemann surface,  and call $\Ggot$ its Gauss map.  Assume that either of the following statements holds:
\begin{itemize}
\item $M$ carries non-constant bounded holomorphic functions, $X$ is asymptotically conformal  and  $\Ggot(M) \subset \s^2\cap\{|x_3|<1-\epsilon\}$ for some $\epsilon>0.$ 
\item  $M=\d,$  $\overline{\c}-g(\d)$ has positive logarithmic capacity and the Beltrami differential $\mu$ of $\Ggot$ satisfies the quadratic Carleson measure condition 
$$ \text{$\lim_{R \to 0} \nu_R=0$ and $\int_{0}^1 \nu_R^2 R^{-2} dR<+\infty,$}$$ where $\nu_R=\sup\{|\mu(x)|\;|\; x \in \Gamma_R\},$ $\Gamma_R$ is the set of points in $\d$ at distance $R>0$ from $\Gamma,$ and $\Gamma \subset \partial(\d)$ is an open arc.
\end{itemize}
Then $X$ is conformal and incomplete.}
\end{quote}

This paper is laid out as follows. In Section \ref{sec:prelim} we introduce some background and notations. We also study the first properties of the Gauss map of harmonic immersions and obtain Weierstrass type formulas for them (see Proposition \ref{pro:wei-h}). Finally, we deal with the basic geometry of QC harmonic immersions. Section \ref{sec:ex} is devoted to present some basic examples of complete harmonic immersions:  graphs over planar domains, rotational surfaces and harmonic horns.  In Section \ref{sec:finite} we study the global geometry of harmonic immersions of FTC and prove Theorem I (see Corollaries \ref{co:fun} and \ref{co:gaussmap} and Theorem \ref{th:anular}). In Section \ref{sec:embe} we classify complete harmonic embedded annuli  with FTC,  construct the examples in Theorem II  and prove their geometrical properties (see Theorems \ref{th:cat-todas}, \ref{th:flujo} and \ref{th:i}). Finally, Theorems III and IV are proved in Section \ref{sec:gauss} (see Theorems \ref{th:mediohuevo}, \ref{th:gaussacotada} and \ref{th:privalov} and Corollary \ref{co:privalov}).


\section{Harmonic immersions} \label{sec:prelim}

Throughout this paper, $\langle,\rangle$ and $\|\cdot\|$ will denote the Euclidean metric  and the Euclidean norm in $\r^3.$ As usual, we write $\c,$ $\overline{\c}=\c \cup \{\infty\}$ and $\d$ for the complex plane, the Riemann sphere and the unit complex disc, respectively. Moreover,  $M$ will denote either an open Riemann surface or a Riemann surface with non-empty compact boundary. 

An open Riemann surface is said to be {\em hyperbolic} if an only if it carries a negative non-constant subharmonic
function. Otherwise, it is said to be {\em parabolic}. If $M$ has non-empty boundary then $M$ is said to be parabolic if bounded
harmonic functions on $M$ are determined by their boundary values, or equivalently, if the harmonic
measure of $M$ respect to a point $P\in M-\partial(M)$ is full on $\partial(M).$ 
$M$ is said to be of {\em finite conformal type} if $M=\overline{M}-\{P_1,\ldots,P_r\},$ where $\overline{M}$ is a compact Riemann surface possibly with non empty boundary $\partial(M)$ and $\{P_1,\ldots,P_r\}\subset \overline{M}-\partial (M).$ In this case $\overline{M}$ is said to be the {\em compactification} of $M.$ 
$M$  is said to be an {\em annular end} if it is homeomorphic to $[0,1[ \times \s^1.$ In this case and up to biholomorphisms,  $M=\{r<|z|\leq 1\},$ where $r \in [0,1[.$ Notice that  $M$ is parabolic (or of finite conformal type) if and only if $r=0,$ and in this case $\overline{M}:=\overline{\d}.$ 
A smooth 1-form $\tau$ on $M$ is said to be of type (1,0) if for any conformal chart (U,z) on $M,$ $\tau(z)=f(z) dz$ for a smooth function $f:U \to \c.$

See \cite{ahlfors-sario} for a good reference on Riemann surfaces theory.

An immersion $X:M \to \r^3$ is said to be {\em conformal} if the metric $X^*(\langle,\rangle)$ is conformal on $M.$

Let us present the fundamental objects of study in this article.
\begin{definition} 
A map $X=(X_j)_{j=1,2,3}:M \to \r^3$ is said to be a {\em harmonic immersion} if $X$ is an immersion and $X_j$ is a harmonic function on $M,$ $j=1,2,3.$ 
A subset $\Scal \subset \r^3$ is said to be a {\em harmonic surface} if there exists a harmonic immersion  $X:M \to \r^3$ such that $\Scal=X(M).$ In this case, $X$ is said to be a harmonic parameterization of $\Scal.$
\end{definition}
A harmonic surface may admit different harmonic parameterizations. For instance, 
\begin{equation} \label{eq:helicoid}
\begin{array}{rll}
Y_1:\c\to\r^3,&&Y_1(z)=\Re(e^z,\imath e^z,\imath z), \;\;\text{and}\\
Y_2:\{\Re(z)>0\} \to \r^3,&&Y_2(z)=\Re\big(\sinh(z), \imath \cosh(z),\imath z\big)\\
\end{array}
\end{equation}
are two different harmonic parameterizations of the same surface (the half helicoid). Notice that $Y_2$ is conformal, whereas $Y_1$ is not.

\begin{remark}
A harmonic surface has no elliptic points by the maximum principle for harmonic functions. In other words,  its Gauss curvature is nowhere positive.
\end{remark}
In the sequel, $X:M \to \r^3$ will denote a harmonic immersion.

Given a linear transformation $A\in Gl(3,\r),$ the map   $A\circ X:M\to\r^3$ is a harmonic immersion as well. Harmonic immersions differing in linear transformations and translations  are said to be {\em linearly equivalent}.

Denote by $\partial_z$ and $\partial_{\overline{z}}$  the global complex operators on $M$ given by $\partial_z|_U=\frac{\partial}{\partial w} dw$ and $\partial_{\overline{z}}|_U=\frac{\partial}{\partial \overline{w}} dw,$ respectively, for any conformal chart $(U,w)$ on $M.$ Put $X=(X_j)_{j=1,2,3}$ and set $$\Phi=(\Phi_j)_{j=1,2,3}:=(\partial_z X_j)_{j=1,2,3}.$$ Notice that $\Phi_j$ is a holomorphic 1-form on $M$ for all $j.$ 
\begin{definition} \label{def:hopf}
The couple $(M,\Phi)$  is called the {\em Weierstrass representation} of $X.$ 
The holomorphic 2-form $\Hgot:=\sum_{j=1}^3 \Phi_j^2$ on $M$ is said to be the {\em Hopf differential} of the harmonic immersion $X.$
\end{definition} 
Let $w=u+\imath v$ be a local conformal parameter on $M$ and write $X_u=\frac{\partial X}{\partial u},$ $X_v=\frac{\partial X}{\partial v}.$ 
Notice that $|\Hgot(w)|^2=\left( (\|X_u\|^2-\|X_v\|^2)^2+4 \langle X_u,X_v\rangle^2 \right) |dw|^4,$ where $\|\cdot \|$ is the Euclidean norm. From Cauchy-Schwarz inequality, 
$|\Hgot(w)| < (\|X_u\|^2+\|X_v\|^2) |dw|^2=\sum_{j=1}^3 |\Phi_j|^2(w),$ where we have taken into account that $\{X_u,X_v\}$ are linearly independent. 

In other words, if we set $\|\Phi\|:=\left(\sum_{j=1}^3 |\Phi_j|^2\right)^{1/2}$  then
\begin{equation} \label{eq:condimmer}
|\Hgot|<\|\Phi\|^2\;\; \mbox {everywhere on}\;\; M
\end{equation}
and $\|\Phi\|^2$ is a conformal metric on $M.$
The converse holds:
\begin{lemma} \label{lem:wei}
Let $\Phi=(\Phi_j)_{j=1,2,3}$ be a triple of holomorphic 1-forms on $M$ such that
$|\sum_{j=1}^3 \Phi_j^2|<\|\Phi\|^2$ everywhere on $M$ and $\Phi$ has no real periods on $M.$

Then the map $X:M \to \r^3,$ $X(P)=\Re \left(\int^P \Phi \right),$ is a harmonic immersion.
\end{lemma}

\begin{definition}
The conformal Riemannian metric $\|\Phi\|^2$ is said to be the {\em Klotz metric} associated to $X.$ The immersion $X$ is said to be {\em weakly complete} if its associated Klotz metric $\|\Phi\|^2$ is complete. 
\end{definition}

\begin{remark} \label{re:kcompleta} By \cite{K}, $X^*(\langle,\,\rangle)\leq \|\Phi\|^2.$ In particular, if $X$ is  complete then  $X$ is weakly complete (the converse is not true, see the example in Remark \ref{re:contra} below). 
\end{remark}

\subsection{Gauss map and curvature.}

Let $\Ggot:M \to \s^2$ be the orientation preserving Gauss map  of $X.$ If $(U,z)$ is a conformal parameter in $M$ and $\Phi_j(z)=\phi_j(z) dz,$ then 
\begin{equation} \label{eq:gauss}
\Ggot|_U=\frac{\Im (\phi_2 \overline{\phi}_3,\phi_3 \overline{\phi}_1,\phi_1 \overline{\phi}_2)}{\| \Im (\phi_2 \overline{\phi}_3,\phi_3 \overline{\phi}_1,\phi_1 \overline{\phi}_2) \|}=\frac{\imath \Phi\wedge \overline{\Phi}}{\|\imath \Phi\wedge \overline{\Phi}\|}.
\end{equation}
Therefore, $\Ggot$ is an analytical map. If we denote by $\xi:\s^2-\{(0,0,1)\} \to \r^2,$ $\xi(x_1,x_2,x_3)=(x_1/(1-x_3),x_2/(1-x_3))$  the stereographic projection, then the orientation preserving map $g:=\xi \circ \Ggot$ is said to be the {\em complex Gauss map} of $X.$ 

The topological properties of $\Ggot$ are described in the following
\begin{lemma} \label{lem:gauss}
Assume that $X(M)$ does not lie in a plane. Then the following assertions hold:
\begin{enumerate}[{\rm (a)}]
\item For any $\nu \in \s^2,$ $\Ggot^{-1}(\{\nu,-\nu\})$ coincides with the zero set of  $\langle \nu,\Phi\rangle$ on $M.$
\item The Gauss map $\Ggot:M \to \s^2$ is open and $d\Ggot$ is bijective except for a discrete subset $S$ of points in $M.$
\item If $P \in S$ then there exists an open disc $U \subset M$ containing $P$ such that $\Ggot|_U:U \to \Ggot(U)$ is a $k$-sheeted covering ramified only at $P,$ where  $k$ is the zero order of $\langle \Ggot(P),\Phi\rangle$ at $P.$
\end{enumerate}
\end{lemma}
\begin{proof} Take $\nu \in \s^2$ and consider a rigid motion $A$ in $\r^3$ preserving the orientation of both $M$ and $\r^3,$ and with $A(\nu)=(0,0,1).$  By equation \eqref{eq:gauss},
$\langle A \circ \Phi,(0,0,1)\rangle(P)=0$ if and only if $A \circ \Ggot(P)\in\{(0,0,1),(0,0,-1)\},$   which proves (a).

Let us check (b) and (c). Take $P \in M,$ and up to a rigid motion as above assume that $\Ggot(P)=(0,0,-1).$ By equation \eqref{eq:gauss} $\Phi_3(P)=0,$ and since $X(M)$ does not lie in a horizontal plane then $\Phi_3$ is not everywhere zero.  Let $(U,z)$ be a conformal disc in $M$ centered at $P,$ and write $\Phi_j(z)=\phi_j(z) dz,$ $j=1,2,3.$ Then we can put
$\phi_1(z)=(c_1+z h_1(z)),$ $\phi_2(z)=(c_2+z h_2(z))$ and $\phi_3(z)=z^k h_3(z),$ where  $c_j\in \c,$ $j=1,2,$  $h_j(z)$ is a holomorphic function, $j=1,2,3,$ $c_3:=h_3(0) \neq 0$ and $k$ is the natural number given in the statement of the lemma.  Notice that $\Im(c_1 \overline{c_2})\neq 0$ (otherwise $X$ would  not be an immersion at $P$). Up to a positive linear isomorphism preserving the $x_3$-axis, we can  assume that $\overline{c}_2 c_3=-\imath$ and $\overline{c}_1 c_3=1$ (recall that $g$ is orientation preserving).

Therefore $g(P)=0,$ and by equation \eqref{eq:gauss}
\[
g(z)=
\frac{1}{2/|c_3|^2+\Ocal_3(z)}\left(\Re(z^k) +\Ocal_1(z)\;,\; \Im(z^k) +\Ocal_2(z)\right),
\]
where $\limsup_{z \to 0} |\Ocal_j(z)/z^{k+1}|<+\infty,$ $j=1,2,$ and $\limsup_{z \to 0} |\Ocal_3(z)/z|<+\infty.$ 

This shows that $g(U)$ is a neighborhood of $0.$ Furthermore,  $g$ is regular at $P$ if and only if $k=1,$ and otherwise $g$ has an (isolated) topological branch point at $P$ of order $k-1\geq 1.$   Since $P$ is an arbitrary point of $M$ we are done.
\end{proof}

Let us compute the curvatures of $X.$ Consider a conformal parameterization $(U,z=u+\imath v)$ in $M,$ call $\Phi_j|_U=\phi_j(z) dz,$ $j=1,2,3,$ $\Phi=\phi(z) dz,$ and note that $X_u:=\partial X/\partial u=\Re(\phi)$ and $X_v:=\partial X/\partial v=-\Im(\phi).$ Following the classical notation, the first and second fundamental form of $X$ in the basis $\{\Re(\phi),-\Im(\phi)\}$ are determined by the coefficients 
$$E=\langle \Re(\phi),\Re(\phi) \rangle, \;\;F=-\langle \Re(\phi),\Im(\phi) \rangle,\;\; G=\langle \Im(\phi),\Im(\phi) \rangle,\;\; \mbox{and}$$
$$e=-g=\langle \Ggot,\Re(\phi')\rangle,\;\; f=-\langle \Ggot,\Im(\phi')\rangle,$$ where we have chosen the normal direction $\Ggot|_U=-\frac{\Re(\phi)\wedge\Im(\phi)}{\| \Re(\phi)\wedge\Im(\phi)\|}=-\frac{\imath \phi \wedge \overline{\phi}}{\|\imath\phi \wedge \overline{\phi}\|}.$

Since the Gauss curvature $\Kcal$ and the mean curvature $H$ of $X$ are given by $\Kcal=\frac{eg-f^2}{EG-F^2}$ and $H=\frac12\frac{Eg-2Ff+Ge}{EG-F^2},$ straightforward computations give that
\begin{equation}\label{eq:curvatura}
\Kcal=-4\frac{|\langle \Ggot,\phi' \rangle|^2}{\|\phi \wedge \overline{\phi}\|^2},\quad H=-2 \frac{\langle \Ggot,\Re(\overline{h} \phi') \rangle}{\|\phi \wedge \overline{\phi}\|^2},
\end{equation}
where $h:=\Hgot/|dz|^2=\sum_{j=1}^3 \phi_j^2.$ Therefore,
\begin{equation} \label{eq:normsigma0}
\|\sigma\|^2=4 H^2-2\Kcal= \frac{8}{\|\phi \wedge \overline{\phi}\|^2} \left(2 \frac{\langle \Ggot,\Re(\overline{h} \phi') \rangle^2}{\|\phi \wedge \overline{\phi}\|^2}
+|\langle \Ggot,\phi' \rangle|^2 \right),
\end{equation}
where $\|\sigma\|$ is the norm of the second fundamental form of $X.$

By \eqref{eq:curvatura}, \eqref{eq:normsigma0} and Cauchy-Schwarz inequality, one also has that
\begin{equation} \label{eq:K-sigma}
0\leq -2\Kcal\leq \|\sigma\|^2 \leq -2 \Kcal \left(2\frac{|\Hgot|^2}{\|\Phi \wedge \overline{\Phi}\|^2}+1\right).
\end{equation}

\begin{remark}
From \eqref{eq:curvatura}, $K=0$ if and only if $\langle d\Ggot,\phi \rangle=0,$ that is to say, $d\Ggot=0$ and $\Ggot$ is constant. In particular, the only flat harmonic surfaces in $\r^3$ are planes.
\end{remark}

To finish this subsection, let us show that the Weierstrass data of $X$ can be rewritten  in terms of its Hopf differential and complex Gauss map $g.$ Recall that, by Lemma \ref{lem:gauss},  $g$  is an orientation preserving open map with a discrete set of topological branch points. 

\begin{proposition}\label{pro:wei-h}
If we set $\lambda:=\frac{\Phi_3(1+|g|^2)}{2 |g|},$ then the 2-form $\lambda^2- \Hgot$ has a globally defined square root $\eta$ on $M$ so that $\Re(\eta\overline{\lambda}) \geq 0.$ Furthermore,
\begin{enumerate}[{\rm ({W}.1)}]
\item $\displaystyle \Phi=\left(\frac{\Re(g) (1-|g|^2)}{|g|(1+|g|^2)}\lambda -\imath \frac{\Im(g)}{|g|} \eta,\frac{\Im(g) (1-|g|^2)}{|g|(1+|g|^2)}\lambda +\imath \frac{\Re(g)}{|g|} \eta,\frac{2 |g|}{1+|g|^2}\lambda\right),$
\item $\displaystyle \frac{\partial_{\overline{z}} g}{\overline{\partial_z g}} = -\frac{g(\lambda-\eta)}{\overline{g}(\lambda+\eta)},$  and
\item $\displaystyle \frac{|\Hgot|}{|\lambda|^2+|\eta|^2}<1.$
\end{enumerate}
As a consequence, 
\begin{equation} \label{eq:klotz-g} 
\|\Phi\|^2=|\lambda|^2+|\eta|^2.
\end{equation}

Conversely, given two smooth 1-forms  $\eta$ and $\lambda$ of type (1,0) on $M$ with $\Re(\eta\overline{\lambda}) \geq 0,$  and an orientation preserving smooth mapping $g:M \to \overline{\c}$ such that 
\begin{itemize} 
\item $\Phi_3:=2 |g| \lambda/(1+|g|^2)$ and $\Hgot:=\lambda^2-\eta^2$ are holomorphic,
\item  {\rm (W.2)} and {\rm (W.3)} hold,  and
\item the vectorial 1-form $\Phi=(\Phi_j)_{j=1,2,3}$  given by {\rm (W.1)} has no real periods,
\end{itemize}
then $\Phi$ is holomorphic and the map $Y:M \to \r^3,$ $Y=\Re(\int \Phi),$ is a harmonic immersion with Hopf differential $\Hgot$ and complex Gauss map $g.$ 
\end{proposition}
\begin{proof} Since
$\Ggot=\frac1{1+|g|^2}(2 \Re(g),2 \Im(g),|g|^2-1),$  then 
\begin{equation} \label{eq:g}
2 \Re(g) \Phi_1+2 \Im(g) \Phi_2+(|g|^2-1) \Phi_3=0. 
\end{equation}

(W.1) easily follows from  \eqref{eq:g} and the definition of $\Hgot$ (see Definition \ref{def:hopf}) for a suitable branch $\eta$ of $\lambda^2- \Hgot.$ Applying $\partial_{\overline{z}}$ to \eqref{eq:g} we deduce (W.2). Since $g$ is orientation preserving one has that $|\frac{\partial_{\overline{z}} g}{\partial_z g}|\leq 1,$ hence $\Re(\eta\overline{\lambda}) \geq 0$  everywhere on $M.$   Equation \eqref{eq:klotz-g} is an elementary consequence of (W.1). Finally, equations \eqref{eq:klotz-g} and \eqref{eq:condimmer} give (W.3).

For the converse, notice that (W.2) implies that $\partial_{\overline{z}} \Phi=0,$ and so $\Phi$ is holomorphic.  Since $\Phi$ has no real periods, the map $Y$ is well defined and harmonic. Equation \eqref{eq:klotz-g} holds as well, and (W.3) implies that $Y$ is an immersion. The remaining properties follow straightforwardly.  
\end{proof}
If $X$ is conformal (i.e., if $\Hgot=0$) then $\eta=\lambda,$ $g$ is meromorphic and the above  is nothing but the Weierstrass representation for minimal surfaces.

\subsection{Quasiconformal harmonic immersions.}
The geometry of harmonic immersions is strongly influenced by the analytical properties of their Gauss map. For instance, a harmonic immersion $X$ is conformal (and so minimal) if and only if its Gauss map $\Ggot$ does, and in this case $g$ is meromorphic. Since quasiconformality generalizes conformality, it is natural to wonder about the global properties of harmonic immersions with quasiconformal Gauss map. 
Recall that the stereographic projection is conformal, hence $\Ggot$ is quasiconformal if and only if $g$ does, that is to say, if and only if $|\mu|<1-\epsilon,$ $\epsilon >0,$ where $\mu:=\frac{\partial_{\overline{z}} g}{\partial_z g}$ is the Beltrami differential of $g.$

\begin{definition}
A harmonic immersion $X:M \to \r^3$ is said to be {\em quasiconformal} (QC for short) if its orientation preserving Gauss map $\Ggot:M \to \s^2$ is  quasiconformal (or equivalently, if $g$ is quasiconformal). 
In this case, $X$ is said to be a harmonic QC parameterization of the harmonic surface $X(M).$
\end{definition}
A harmonic surface in $\r^3$ may admit both QC and non-QC harmonic parameterizations. Indeed, the immersion $Y_2$  in \eqref{eq:helicoid} is even conformal, whereas $Y_1$ is not QC (for a proof of this assertion, use Proposition \ref{pro:index} below).

\begin{remark}\label{re:beltrami} By Proposition \ref{pro:wei-h}-(W.2), the Beltrami differential of $g$ is given by
\begin{equation} \label{eq:beltrami}
\mu=-\frac{g(\lambda-\eta)}{\overline{g}(\lambda+\eta)} \left(\frac{\overline{\partial_z g}}{\partial_{{z}} g}\right).
\end{equation}
The condition $\Re(\eta\overline{\lambda}) \geq 0$ simply means that $|\mu|\leq 1.$
Therefore $X$ is QC if and only if $\sup_M |\mu|<1.$ 
 
Since $\|X_u\|^2+\|X_v\|^2=|\lambda|^2+|\eta|^2$ and $\|X_u\times X_v\|=\sqrt{\|\Phi\|^2-|\Hgot|^2}=\Re (\overline{\lambda}\eta)$ for any conformal parameter $z=u+\imath v$ on $M,$ then the distortion function of $X$ is given by
\[
\mathfrak{D}_X=\frac{|\lambda|^2+|\eta|^2}{2\Re(\overline{\lambda}\eta}).
\]
Then, $\sup_M|\mu|<1$ if and only if $\mathfrak{D}_X$ is bounded, or in other words, $X$ is QC if and only if $X$ is a quasiconformal map in the classical sense. See \cite{Ka} for an alternative proof of this assertion.
\end{remark}

In Proposition \ref{pro:index} below we will state several equivalent formulations of quasiconformality for harmonic immersions. The following notations are required.

For any bounded function $f:M \to \r,$  write $\limsup_{P \to \infty} f(P)$ for $\lim_{n \in \n} \sup_{M-C_n} f,$ where $C_1 \subset C_2 \subset ...$ is any exhaustion of $M$ by compact sets.  Set 
$$i_X:=\limsup_{P \to \infty} |\mu|\quad \text{and} \quad i^X:=\limsup_{P \to \infty} \frac{|\Hgot|}{\|\Phi\|^2}.$$

\begin{proposition} \label{pro:index}
The inequalities $0\leq i_X \leq i^X\leq \frac{2 i_X}{1+i_X^2}\leq 1$ hold. 

As a consequence, $X$ is QC if and only if either of the following conditions is satisfied:
\begin{enumerate}[\rm (i)]
\item $i^X<1.$
\item $\sup_{M} |\Hgot|/\|\Phi\|^2<1.$
\item $i_X<1.$
\item $\sup_{M} |\mu|<1.$
\end{enumerate}
\end{proposition}
\begin{proof} If we write $y=\frac{\lambda-\eta}{\lambda+\eta},$ equation \eqref{eq:beltrami} gives that $|\mu|=|y|$ and equation \eqref{eq:klotz-g} that $|\Hgot|/\|\Phi\|^2\leq 2 |y|/(1+|y|^2).$ Taking into account that $|y|\leq 1$ (recall that $|\mu|\leq 1$), we deduce  
\begin{equation} \label{eq:beltramis}
|\mu|\leq \frac{|\Hgot|}{\|\Phi\|^2}\leq \frac{2 |\mu|}{1+|\mu|^2}\leq 1\quad \text{on $M.$}
\end{equation}
The first part of the proposition follows from the fact that  $x\mapsto 2x/(1+x^2)$ is increasing in $[0,1].$

 Obviously $X$ is QC if and only if  (iv) holds. Equation \eqref{eq:condimmer} says that (i) $\Leftrightarrow$ (ii). By equation \eqref{eq:beltramis} and the first part of the proposition, (ii) $\Leftrightarrow$ (iv) and  (i) $\Leftrightarrow$ (iii). 
\end{proof}

Let $M_0$ denote the underlying oriented differentiable structure of $M$ and  write  $M=(M_0,\Ccal)$ for a suitable holomorphic atlas $\Ccal$ on $M_0.$ 
Label $\Ecal$ as the holomorphic structure on $M_0$ of isothermal coordinates for the metric $X^*(\langle, \rangle),$ and denote by $M^e$  the Riemann surface $(M_0, \Ecal).$
The real numbers $i_X$ and $i^X$  measure  how far the harmonic immersion is  from being  conformal. In other words and roughly speaking, how far is the conformal structure $\Ccal$  from $\Ecal$ in the moduli space of Riemann surfaces with underlying topology $M_0$ (see Lemma \ref{lem:quasiconformal}-(i) below). In this sense, the statements ${\rm Id}:M \to M^e$ is conformal, $X:M \to \r^3$ is conformal, and $\Hgot=0$ on $M,$ are equivalent.

The notion of QC harmonic immersion is  related to the one of quasiminimal surface introduced by Osserman in \cite{O1}. Recall that a {\em quasiminimal surface} (a QM surface for short) is a $C^2$ surface in $\r^3$ with the property that the principal curvatures  $k_1,$  $k_2$ satisfy an inequality
$
0<\delta\leq-k_1/k_2\leq1/\delta
$
at every point where they do not vanish simultaneously. By the classical formula of Rodrigues, the differential of the Gauss map takes the unit circle into an ellipse
whose semi-minor and semi-major axes are $|k_1|$ and $|k_2|.$ Thus, quasiminimality is equivalent to the assumption that the Gauss curvature is nowhere positive and the Gauss spherical map is quasiconformal with respect to the conformal structure on the surface induced by isothermal charts for the Euclidean metric.

The following lemma provides a new interpretation of quasiconformality, and connects this concept with the one of quasiminimality. 
\begin{lemma}\label{lem:quasiconformal}
The following statements hold:
\begin{enumerate}[\rm (i)]
\item $X$ is QC $\Leftrightarrow$ ${\rm Id}:M \to M^{e}$ is quasiconformal. 
\item If $X$ is QC then $X(M)$ is QM.
\item If $X$ is QC then $-2 \Kcal\leq \|\sigma\|^2 \leq -c \Kcal$ for some constant $c\geq 2.$ In particular, if $X$ is QC then $\int_M |\Kcal| dS<+\infty$ $\Leftrightarrow$ $\int_M \|\sigma\|^2 dS<+\infty,$ where $dS$ is the area element of $X^*(\langle,\rangle).$
\end{enumerate}
\end{lemma}
\begin{proof}  
Let $(U,z=u+\imath v)$ and $(U,\xi=x+\imath y)$ be two conformal charts in $M$ and $M^e,$ respectively, and write $(x(u,v),y(u,v))$ for the transition map.
For a map $f(u,v)$ write $f_u$ and $f_v$  for $\partial f/\partial u$ and $\partial f/\partial v,$ respectively. Label $C_u=(x_u,y_u)$ and $C_v=(x_v,y_v) \in \r^2.$ 

One has that $\|X_u\|^2=\delta^2 \|C_u\|^2,\;\; \|X_v\|^2=\delta^2 \|C_v\|^2$ and $\langle X_u,X_v\rangle=\delta^2\langle C_u,C_v\rangle,$ where we have taken into account that $\xi$ is isothermal and written $\delta=\|X_x\|=\|X_y\|.$ Then $$|\Hgot|^2/\|\Phi\|^4=1-4 \frac{1-\langle C_u/\|C_u\|,C_v/\|C_v\| \rangle^2}{\left(\|C_u\|/\|C_v\|+ \|C_v\|/\|C_u\| \right)^2},$$ and so 
$\sup_{M} \frac{|\Hgot|}{\|\Phi\|^2}<1$ if and only if there exists $\epsilon>0$ such that $ \frac{1-\langle C_u/\|C_u\|,C_v/\|C_v\| \rangle^2}{\left(\|C_u\|/\|C_v\|+ \|C_v\|/\|C_u\| \right)^2}>\epsilon$ for any charts as above, that is to say, if and only if ${\rm Id}:M \to M^{e}$ is quasiconformal. Taking into account Proposition \ref{pro:index}, (i) holds.

To check (ii),  use (i) and notice that if  ${\rm Id}:M \to M^{e}$ is quasiconformal then   $\Ggot:M^e \to \s^2$ is quasiconformal as well.  

To see (iii), note that (ii) yields that $X(M)$ is QM. Thus, $-2 \Kcal\leq \|\sigma\|^2 \leq -c \Kcal$ for some constant $c\geq 2,$  and we are done.  
\end{proof}
\begin{remark}
The converse of item {\rm (ii)} in  Lemma \ref{lem:quasiconformal} is false. For instance, the harmonic immersion $Y_1$ in \eqref{eq:helicoid}  is QM (in fact $Y_1(\c)$ is a piece of the minimal helicoid) but not QC. Indeed, since the Weierstrass data of $Y_1$ are given by $\Phi=(e^z,\imath e^z,\imath)dz,$ then $\frac{|\Hgot|}{\|\Phi\|^2}=\frac{1}{2|e^{2z}|+1},$ and so $i^{Y_1}=1.$
\end{remark}

\begin{definition}
The immersion $X$ is said to be {\em  asymptotically conformal} (AC  for short) if  $i^X=i_X=0.$ Harmonic immersions linearly equivalent to conformal immersions are said to be {\em essentially conformal} (EC for short). Likewise, harmonic immersions linearly equivalent to an AC harmonic immersion are said to be {\em essentially asymptotically conformal} (EAC  for short).
\end{definition}
  
\begin{proposition} \label{pro:linear}
The following statements hold:
\begin{enumerate}[{\rm (a)}]
\item $X$ is conformal $\Rightarrow$ $X$ is AC $\Rightarrow$  $X$ is QC.
\item $X$ is QC $\Leftrightarrow$ $A \circ X$ is QC for all $A \in Gl(3,\r).$
\end{enumerate}
\end{proposition}
\begin{proof} Item (a) is trivial.

To check (b),  take into account Proposition \ref{pro:index} and notice that
$\sup_{M} \frac{|\Hgot|}{\|\Phi\|^2}<1$ if and only if there exists $\epsilon>0$ such that 
\begin{equation} \label{eq:ladelseno}
\epsilon<\|X_u\|/\|X_v\|<1/\epsilon\;\; \text{and}\;\; |\langle X_u/\|X_u\|,X_v/\|X_v\| \rangle|<1-\epsilon
\end{equation} for any conformal chart $(U,z=u+\imath v)$ in $M$ (here $X_u=\partial X/\partial u$ and likewise for $X_v$). For each $A \in Gl(3,\r),$ these inequalities still hold when replacing $X$ and $\epsilon$ for $A \circ X$ and a suitable positive constant $\epsilon_A$ depending on $A$ and $\epsilon,$  respectively. 
\end{proof}
In general, linear transformations preserve neither conformality nor asymptotic conformality. Only orthogonal transformations do. 

\section{Basic examples}\label{sec:ex}

In this section we present some basic and further examples of complete harmonic immersions in $\r^3.$ 
Recall that  any minimal surface in $\r^3$ admits a conformal harmonic parameterization. Although the family of minimal surfaces is very vast and provides a good source of examples, we are mostly interested in the harmonic non-conformal case. 

\subsection{Harmonic graphs}
For any harmonic function $u:\c\to\r,$ the map $X:\c\to\r^3,$ $X(z)=(\Re(z),\Im(z),u(z)),$ is a harmonic immersion  and $X(\c)$ is an entire graph over the plane $\{x_3=0\}.$ The plane is the simplest case (choose $u=0$).

With a little more effort we can construct complete harmonic graphs over the unit disc.  
We will need the following technical lemma.
\begin{lemma}\label{lem:JX}
There exists a harmonic function $u:\d\to\r$ such that $\int_\alpha |\Re(\partial_z u)|=+\infty$  for any divergent curve $\alpha \subset \d$ with finite Euclidean length.
\end{lemma}
\begin{proof}
Consider the labyrinth of compact sets $\{K_n\}_{n\in\n}$  in $\d$ described by Jorge-Xavier in \cite{j-x}. Let $r_n>1$ be the quotient between the outer and the inner radii of the truncated annulus $K_n,$ $n\in\n,$ and notice that $\lim_{n \to \infty} r_n=1.$  Let $h:\d\to\c$ be a holomorphic function satisfying 
\begin{equation}\label{eq:laberinto}
\left|h'(z)-\frac{1}{\ln(r_n) z}\right|<1,\quad \forall z\in K_n,\quad \forall n\in\n.
\end{equation}
The existence of such a function $h$ is implied by Lemma 2 in \cite{j-x}.

Let $\alpha \subset \d$ be a divergent curve with finite Euclidean length. Without loss of generality we can assume that $\alpha$ crosses all the compact sets $\{K_n\;|\; n$ even, $n\geq n_0\},$ for some $n_0\in\n$ (see \cite{j-x} for details). From \eqref{eq:laberinto}, $\forall n$ even, $n\geq n_0,$
\[
\int_{\alpha\cap K_n} |\Re(h'(z)dz)|\geq \left| \int_{\alpha\cap K_n}\left( \left| \Re\left( \frac{dz}{\ln(r_n) z}\right)\right|-|dz|\right)\right| \geq |1-\text{length} (\alpha\cap K_n)|.
\]

This inequality and the finiteness of the Euclidean length of $\alpha$ give
\[
\int_{\alpha} |\Re(h'(z)dz)|\geq \left(\sum_{n\text{ even, }n\geq n_0} 1\right)-\text{length}(\alpha)=+\infty.
\]
To finish, it suffices to choose $u=\Re(h).$
\end{proof}

\begin{corollary}
There exist harmonic functions $u:\d\to\r$ so that the map $X:\d\to\r^3,$ $X(z)=(\Re(z),\Im(z),u(z)),$ is a complete harmonic immersion.
\end{corollary}
\begin{proof} Let $u:\d \to \r$ be a harmonic function given by Lemma \ref{lem:JX}, and set $X:\d\to\r^3$ as in the statement. Let $\alpha \subset \d$ be a divergent curve. If $\alpha$ has infinite Euclidean length then the same holds for $X(\alpha)$ (recall that $X(\d)$ is a graph).  Lemma \ref{lem:JX} implies that $X(\alpha)$ has infinite length when $\alpha$ has finite Euclidean length.  
\end{proof}
Complete harmonic graphs over finitely connected planar domains can be constructed by similar arguments.

\subsection{Rotational harmonic surfaces}  Let $X=(X_j)_{j=1,2,3}:M\to\r^3$ be a non-flat rotational harmonic annulus. Without loss of generality, assume that $X(M)$ is invariant under the family of vertical rotations $\{B_\theta:\r^3\to\r^3\}_{\theta\in\r},$ $B_\theta(x,y,t)=(\cos(\theta)x+\sin(\theta)y, -\sin(\theta)x+\cos(\theta)y,t).$  Write  $(\Phi_j)_{j=1,2,3}=(\partial_z X_j)_{j=1,2,3},$ and let $m\neq 0$ denote the additive period $\imath \Im (\int_\gamma \Phi_3)$ of the  conjugate harmonic function $X_3^*:M \to \r$ of $X_3$ (here $\gamma$ is the loop generating $H_1(M,\z)$).  Set $z=e^{(X_3+\imath X_3^*)\frac{2 \pi}{m}},$ and observe that up to natural identifications  $M=\{r_1<|z|<r_2\},$ where $0<r_1<r_2<+\infty,$ and $\Phi_3(z)=dz/z.$   We also write $\Phi_j=\phi_j(z)dz,$ $j=1,2.$ Up to a symmetry with respect to a vertical plane, we can suppose that $B_\theta$ induces on $M$ the biholomorphism $z\mapsto e^{\imath \theta}z.$ Then, one has
\[
\left.
\begin{array}{r}
\cos(\theta)\phi_1(z)+\sin(\theta)\phi_2(z)=e^{\imath\theta}\phi_1(e^{\imath\theta}z)\\
-\sin(\theta)\phi_1(z)+\cos(\theta)\phi_2(z)=e^{\imath\theta}\phi_2(e^{\imath\theta}z)
\end{array}\right\} \quad \forall \theta\in\r,\; z\in M,
\]
and so $\phi_1(e^{\imath\theta}z)-\imath\phi_2(e^{\imath\theta}z)=\phi_1(z)-\imath\phi_2(z)$ and $\phi_1(e^{\imath\theta}z)+\imath\phi_2(e^{\imath\theta}z)=e^{-2\imath\theta}(\phi_1(z)+\imath\phi_2(z)),$ $\forall \theta\in\r,$ $z\in M.$
This gives that $\phi_1-\imath \phi_2$ and $z^2(\phi_1+\imath \phi_2)$ are constant. Up to the change $z \to c z$ or  $z \to c/z$ for a suitable $c \in \c-\{0\},$ one has
\begin{equation}\label{eq:rota}
\Phi_1=\left(1-\frac{b}{z^2}\right)dz,\quad \Phi_2=\imath\left(1+\frac{b}{z^2}\right)dz,\quad \Phi_3=\frac{1}{z}dz,\quad \Hgot=\frac{1-4b}{z^2}dz^2,
\end{equation}
for a constant $b\in\c.$

\begin{proposition} \label{pro:ab}
$X$ is a complete immersion if and only if $M=\c-\{0\}$ and  $b \notin (-\infty,0) \subset \r.$
\end{proposition}
\begin{proof}
If $X$ were a complete immersion,  \eqref{eq:rota} would imply that $M=\c-\{0\}.$ 
On the other hand, Lemma \ref{lem:wei} and \eqref{eq:rota} give that $X$ is immersion if and only if
\begin{equation}\label{eq:immer}
2 |z|^4 +(1- |1 - 4 b|) |z|^2 + 2 |b|^2>0,\quad\forall z\in\c-\{0\}.
\end{equation}
If $b=0,$ the  inequality \eqref{eq:immer} holds (and so $X$ is an immersion). When $b\neq 0,$ inequality \eqref{eq:immer} holds if and only if the function $p:\r\to\r,$ $p(x)=2 x^4 +(1- |1 - 4 b|) x^2 + 2 |b|^2,$ never vanishes. This happens if and only if either $|1-4b|\leq 1$ or $|1-|1-4b||<4|b|,$ or equivalently $b \notin \{z \in \c\; |\; \Re(z)<0,\; \Im(z)=0\}.$
\end{proof}
Denote by $X_{b}:\c-\{0\}\to \r^3$  the harmonic immersion with Weierstrass representation given by \eqref{eq:rota} for $b$ as in Proposition \ref{pro:ab}.
The classification of complete rotational harmonic surfaces in $\r^3$ can be found in the following
\begin{corollary}
Let $X:M \to \r^3$ be a  non-flat complete harmonic immersion, and assume that $X(M)$ is a rotational surface. 

Then, up to coverings, homotheties and rigid motions,  $X=X_{b}$  for some  $b \notin (-\infty,0).$ 
\end{corollary}
$X_0$ (respectively, $X_b,$ $b\neq 0$) is said to be the {\em rotational harmonic horn} (respectively, a {\em rotational harmonic catenoid}), see Figure \ref{fig:chincheta-catenoide}.

\begin{remark} Elementary computations show that:
\begin{itemize}
\item The rotational harmonic horn  is not QC.
\item Any rotational harmonic catenoid is AC. Furthermore,
\begin{itemize}
\item  two rotational harmonic catenoids $X_{b_1}$ and $X_{b_2}$ are linearly equivalent $\Leftrightarrow$  $b_1=b_2.$
\item  $X_{b}$ is EC $\Leftrightarrow$ $b\in\r$ and $b>0.$
\item  $X_{b}(\c-\{0\})$ is a minimal surface $\Leftrightarrow$ $b=1/4$ $\Leftrightarrow$ $X_{b}$ is conformal.
\end{itemize}
\end{itemize}
\end{remark}
The rotational harmonic horn lies in the boundary of the moduli space of rotational harmonic catenoids. The last proposition shows that in fact it is limit of EC rotational catenoids. The  family $\{X_{b}\;|\; b\notin\r\}$ consists of "purely" harmonic QC immersions, in the sense that they are  linearly equivalent to neither  conformal nor  minimal immersions.

\begin{figure}[ht]
    \begin{center}
    \scalebox{0.25}{\includegraphics{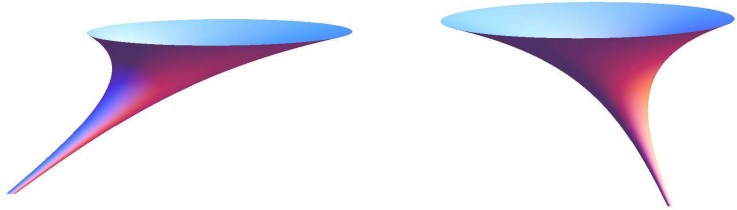}}
        \end{center}
\caption{Two harmonic horns with non-vertical flux.}\label{fig:general-horns}
\end{figure}

The rotational harmonic horn lies in a larger family of complete embedded harmonic annuli. By definition, a  {\em harmonic horn} is any of the harmonic immersions in the following family 
$$X:\c-\{0\} \to \r^3,\quad X(z):=\big(r_1 \log(|z|)+\Re(z), r_2 \log(|z|)+\Im(z),\log(|z|) \big),$$ where $r_1,$ $r_2 \in \r.$ When $(r_1,r_2) \neq (0,0),$ the corresponding harmonic horn has non-vertical flux (see Figure \ref{fig:general-horns}).

\section{Harmonic immersions of finite total curvature} \label{sec:finite}

In this section we will study the global geometry of harmonic immersions of finite total curvature. Our main results are Theorem \ref{th:fun} and Corollary \ref{co:fun}, where this objects are characterized from different points of view. Most of the subsequent corollaries  are extensions of some well-known results for minimal surfaces in $\r^3.$

\begin{definition}
A harmonic immersion $X:M \to \r^3$ is said to be {\em algebraic} if  $M$ has finite conformal type and $\Phi$ extends meromorphically to  $\overline{M}$  with poles at all the ends of $M.$ 

On the other hand, $X$ is said to be of {\em finite total curvature} (FTC  for short) if $\int_M \|\sigma\|^2 dS<+\infty,$ where $\sigma$ is the second fundamental form of $X$ and $dS$ is the area element of $X^*(\langle,\rangle).$
\end{definition}
Obviously, $\int_M |\Kcal| dS<+\infty$ is a much weaker condition than FTC.
\begin{remark} \label{re:contra}
There exist {\em non-complete} algebraic harmonic immersions. For instance, take $X:\c-\{0\} \to \r^3,$ $X(z)=\Re\left(z,\imath z,1/z\right),$ and  notice that the divergent curve $c:[-1,0[ \to \c-\{0\},$ $c(t)=\imath t,$ has finite length in $(\c-\{0\},X^*(\langle,\rangle)).$
\end{remark}

\begin{definition} \label{def:iq} Let $X:M \to \r^3$ be an algebraic harmonic immersion, and write $\overline{M}$ and $\Phi=(\Phi_j)_{j=1,2,3}$ for the compactification of $M$ and the Weierstrass data of $X,$ respectively. For each $P \in \overline{M}-M,$ we set
$I_P:=\max\{{\rm Ord}_P (\Phi_j)\;|\; j=1,2,3\}-1\geq 0,$ where ${\rm Ord}_P(\cdot)$ means pole order at $P.$ The natural number $I_P$ is said to be the {\em weight} of the end $P.$
\end{definition}

If $X$ is a complete harmonic immersion with FTC of an open Riemann surface $M$ into $\r^3,$ then Huber's theorem \cite{huber} says that $M$ has finite topology, and the conformal structure on $M$ induced by the (positively oriented) isothermal charts for $X$ (which in general is different from the original one of $M$) is parabolic. Therefore,  this classical theorem gives no information about the conformal type of $M$ itself. This question and more are answered in the following theorem and its corollaries.

\begin{theorem} \label{th:fun}
Let $A$ be a conformal annular end, let   $X=(X_j)_{j=1,2,3}:A \to \r^3$ be a harmonic immersion, and let $\Phi=(\Phi_j)_{j=1,2,3}$ denote its Weierstrass data. Then the following statements are equivalent:
\begin{enumerate}[\rm (i)]
\item $X$ is complete and of FTC.
\item $X$ is complete and algebraic, and the Gauss map $\Ggot:A \to \s^2$ of $X$ extends continuously to the compactification $\overline{A}\equiv \overline{\d}$ of $A.$
\item  $X$ is algebraic, and if $Q_0\in \overline{A}$ is the unique end of $A,$ then   $I_{Q_0} \geq 1$ and there exists  $R \in O(3,\r)$ such  that 
$$\text{${\rm Ord}_{Q_0}\left((R \circ \Phi)_3\right)<{\rm Ord}_{Q_0}\left((R \circ \Phi)_2\right)={\rm Ord}_{Q_0}\left((R \circ \Phi)_1\right)=I_{Q_0}+1$ and $\frac{(R \circ \Phi)_2}{(R \circ \Phi)_1}({Q_0}) \notin \r.$}$$  
\item $X$ is complete, algebraic and QC.
\item $X$ is complete, QC and  $\int_A |\Kcal| dS<+\infty$.   
\end{enumerate}
Furthermore, if the above statements hold then $X$ is proper, $\|\sigma\|^2$ is bounded,  $X$ is EAC and, up to removing a compact subset, $X(A)$ is a $I_{Q_0}-$sheeted sublinear multigraph over $\Pi-D,$ where $\Pi$ is the limit tangent plane of $X$ at $Q_0$ and $D \subset \Pi$ is an open disc. 
\end{theorem}
\begin{proof} 
Let us see (i) $\Rightarrow$ (ii). Let $A_0$ denote the underlying {\em differentiable} structure of $A$ and  write  $A=(A_0,\Ccal)$ for a suitable holomorphic atlas $\Ccal$ on $A_0.$ 
Likewise, call $\Ecal$ the conformal structure on $A_0$ of isothermal coordinates for the metric $X^*(\langle, \rangle),$ and label $A^e$ as the Riemann surface $(A_0, \Ecal).$ 

Taking into account \eqref{eq:K-sigma} and combining results by Huber \cite{huber} and White \cite{W}, we infer that  $X$ is proper, $A^e$ is biholomorphic to $\overline{\d}-\{0\}$ and the Gauss map $\Ggot:A^e\to\s^2$ of  $X$ extends continuously to all of $\overline{A}^e \equiv \overline{\d}.$ 

Label $P_0\in\overline{A}^e$ as the unique end of $A^e,$ and up to a rigid motion assume that $\Ggot(P_0)=(0,0,1).$ Without loss of generality, suppose that $\Ggot(A^e) \subset  \s^2 \cap \{x_3\geq 1/2\}.$ Take $R>0$ large enough so that $X(\partial(A^e)) \subset C_R:=\{(x_1,x_2,x_3)\,|\;x_ 1^2+x_2^2<R\}.$ By Jorge-Meeks \cite[Theorem 1]{JM},    $A^e \cap X^{-1}(\overline{C}_R)$ is connected and compact,  $A':=A^e-X^{-1}(C_R)$ is connected and  
\begin{equation} \label{eq:sheet}
\text{$(X_1,X_2)|_{A'}:A' \to \{(x_1,x_2)\,|\;x_ 1^2+x_2^2\geq R\}$ is a finite covering of sublinear growth.}
\end{equation}  
In the sequel we relabel $A^e=A'.$ As a consequence,  $X^{-1}(\{x_j=t\})$ consists of a finite number of regular curves $\forall t \in \r,$ $j=1,2.$ If $A$ were hyperbolic, then $(X_1,X_2)^{-1}([0,+\infty)^2)$ would contain a hyperbolic connected component $U$ (see \cite{MP}), contradicting that $X_1+X_2:U \to \r$ is a positive proper harmonic function.  

Therefore,  $A$ is parabolic as well, and using  Meeks and P\'{e}rez results \cite{MP} again,   $\Phi_j$ extends meromorphically to the compactification $\overline{A}\equiv \overline{\d}$ of $A,$ $j=1,2.$ To finish, let us see that $\Phi_3$ extends meromorphically to $\overline{A}.$  Indeed, otherwise $\Phi_3$ would have an essential singularity at the unique end $Q_0$ of $A,$ hence $\Phi_3/\Phi_j$ would have an essential singularity at $Q_0$ as well, $j=1,2.$ By the great Picard theorem, there would exist a sequence $\{P_n\}_{n \in \n} \to Q_0$ such that $\lim_{n \to \infty} \Ggot(P_n)$ was horizontal (see equation \eqref{eq:gauss}), contradicting that $\Ggot(Q_0)=(0,0,1).$

This shows that $X$ is algebraic and (ii) holds. 

Let us check that (ii) $\Rightarrow$ (iii). Let $R$ be a linear isometry in $\r^3$ preserving the orientation of both $A$ and $\r^3,$ and such that $(R \circ \Ggot)(Q_0)=(0,0,1).$ Write $Y=R \circ X,$ $\Psi=(\Psi_j)_{j=1,2,3}=R\circ \Phi$ and $G=R \circ \Ggot,$ and notice that $Y$ is a complete harmonic immersion satisfying  (ii) with $G(Q_0)=(0,0,1).$  

Assume for a moment that ${\rm Ord}_{Q_0}(\Psi_3)\geq {\rm Ord}_{Q_0}(\Psi_j)$ for some $j\in \{1,2\}$ (without loss of generality, for $j=1$).  Then it is not hard to check that 
$$\limsup_{Q \to {Q_0}} \frac{|\Im(\Psi_2 \overline{\Psi}_3)|}{|\Im(\Psi_1 \overline{\Psi}_2)|}>0,$$ hence  equation \eqref{eq:gauss} gives that $G({Q_0})\neq (0,0,1),$ a contradiction.

Let us check now that  ${\rm Ord}_{Q_0}(\Psi_1)={\rm Ord}_{Q_0}(\Psi_2)$ and $(\Psi_2/\Psi_1)({Q_0})\notin \r.$ Otherwise and up to a linear transformation in $\r^3$ fixing pointwise the $x_3$-axis, we can suppose that  ${\rm Ord}_{Q_0}(\Psi_1)>{\rm Ord}_{Q_0}(\Psi_2)>{\rm Ord}_{Q_0}(\Psi_3).$ Then, we can find a sequence $\{P_n\}_{n \in \n}$ in $A$ converging to ${Q_0}$ such that $\Im(\Psi_1\overline{\Psi}_2)(P_n)=0$ for all $n \in \n.$ This shows that $G({Q_0})\neq (0,0,1)$ by  equation \eqref{eq:gauss}, which is absurd.

To finish, let us show that $I_{Q_0}>0.$ Indeed, if ${\rm Ord}_{Q_0}(\Psi_1)={\rm Ord}_{Q_0}(\Psi_2)=1$ and since $(\Psi_2/\Psi_1)({Q_0})\notin \r,$ then either $\Psi_1$ or $\Psi_2$ has non purely real residue at ${Q_0}.$ Thus, $Y$ would not be well defined and (iii) holds. Notice that (iii) and \eqref{eq:sheet} imply that, up to removing a compact subset,  $X(A)$ is a $I_{Q_0}$-sheeted sublinear multigraph over the complement of a bounded set of $\{x_3=0\}.$

Let us see that (iii) $\Rightarrow$ (iv). The  algebraic conditions  in (iii) clearly imply that $X$ is complete. Then it suffices to prove that $X$ is QC.  Fix a conformal parameter $z:A \to \overline{\d}-\{0\}$ and identify   $A=\overline{\d}-\{0\}.$ By item (iii), there exists a linear transformation $L$ preserving the $x_3$-axis so that the Weierstrass data $\Psi=(\Psi_j)_{j=1,2,3}$ of $Y:=L \circ X$ are given by 
\[
\Psi_1(z)=\frac{1+\Ocal_1(z)}{z^{n+m}}dz,\quad \Psi_2(z)=\frac{\imath+\Ocal_2(z)}{z^{n+m}}dz\quad\text{and}\quad \Psi_3(z)=\frac{c+\Ocal_3(z)}{z^n}dz,
\]
where $n,$ $m\in\n,$  $c\in\c-\{0\}$ and $\lim \sup_{z\to 0}|\Ocal_j(z)/z|<+\infty$ $\forall j=1,2,3.$ 

Then $\limsup_{z\to 0} \frac{|\Hgot_Y|}{\|\Psi\|^2}(z)=0,$  where $\Hgot_Y$ is the Hopf differential of $Y.$ Thus  $Y$ is AC, hence $X$ is EAC and so QC (see Proposition \ref{pro:linear}) and we are done.

Let us check now that (iv) $\Rightarrow$ (ii). 
Label $\Hgot$ and $g$ as the Hopf differential and the complex Gauss map of $X,$ respectively. Up to a suitable linear transformation  (which preserves both algebraicity and  quasiconformality, see Proposition \ref{pro:linear}), we can assume that $\lim_{Q \to {Q_0}}\Hgot/\Phi_3^2(Q)=\infty.$ As $X$ is QC then there exist a quasiconformal homeomorphism $q:\overline{A} \to \overline{A}$ fixing ${Q_0}$ and a meromorphic function $f:A \to \c$ such that $g=f \circ q.$ Reason by contradiction and assume that $g$ does not extend continuously to the unique end ${Q_0}$ of $A.$ Then $f$ has an essential singularity at ${Q_0},$ and by the big Picard Theorem, there exists a sequence $\{P_n\}_{n \in \n} \to {Q_0}$ such that $\{|f(P_n)|+1/|f(P_n)|\}_{n\in \n}$ (and so, $\{|g(q^{-1}(P_n))|+1/|g(q^{-1}(P_n))|\}_{n\in \n}$) is bounded.  This implies that 
$\{(\Hgot/\lambda^2) (q^{-1}(P_n))\}_{n \in \n} \to \infty,$ where $2\lambda=\Phi_3(|g|+1/|g|),$ hence the immersion $X$ is not QC  by  \eqref{eq:beltrami}, a contradiction.

Let us check that (iv) $\Rightarrow$ (v). We only have to prove that  $\int_A |\Kcal| dS<+\infty.$ As above, take a conformal parameter $z:A \to \overline{\d}-\{0\}$ and identify   $A=\overline{\d}-\{0\}.$ Write $\Phi(z)=\phi(z) dz.$ 
By equation \eqref{eq:curvatura} and Remark \ref{re:kcompleta}, one has that
$$\int_A |\Kcal| dS= \int_{A} 4\frac{|\langle \Ggot,\phi' \rangle|^2}{\|\phi \wedge \overline{\phi}\|^2} dS \leq \int_{A} 4\frac{|\langle \Ggot,\phi' \rangle|^2}{\|\phi \wedge \overline{\phi}\|^2} \|\phi\|^2 |dz|^2.$$ Since $|\Hgot|^2=\|\Phi\|^4-\|\Phi\wedge \overline{\Phi}\|^2$  and $\sup_{A} \frac{|\Hgot|}{\|\Phi\|^2}<1$ (see Proposition \ref{pro:index}), the metrics $\|\Phi \wedge \overline{\Phi}\|$ and $\|\Phi\|^2$ are equivalent. Thus,  
$$\int_A |\Kcal| dS < a  \int_{A} \frac{|\langle \Ggot,\phi' \rangle|^2}{\| \phi\|^2}|dz|^2$$ for a suitable constant $a>0.$
Then, it suffices to prove that $\int_{A} \frac{|\langle \Ggot,\phi' \rangle|^2}{\| \phi\|^2}|dz|^2<+\infty.$ Taking into account that (iv) $\Rightarrow$ (ii) $\Rightarrow$ (iii), and up to a suitable linear transformation (which preserves the finiteness of the total curvature), we can suppose  that 
\begin{equation}\label{eq:datos}
\Phi_1(z)=\frac{1+\Ocal_1(z)}{z^{n+m}}dz,\quad \Phi_2(z)=\frac{\imath+\Ocal_2(z)}{z^{n+m}}dz\quad\text{and}\quad \Phi_3(z)=\frac{c+\Ocal_3(z)}{z^n}dz,
\end{equation}
where $n,$ $m\in\n,$ $c\in\c-\{0\}$ and $\lim \sup_{z\to 0}|\Ocal_j(z)/z|<+\infty$ $\forall j=1,2,3.$ From \eqref{eq:datos} and equation \eqref{eq:gauss}, one has that $1/\|\phi(z)\|^2\leq c_1 |z|^{2(n+m)}$ and $|\langle \Ggot(z),\phi'(z) \rangle|^2\leq c_2/|z|^{2(n+1)}$ for positive constants $c_1$ and $c_2,$ proving that $\Kcal$ is bounded and (v). By Lemma \ref{lem:quasiconformal}-(iii), $\|\sigma\|^2$ is bounded as well.

Finally, (v) $\Rightarrow$ (i) follows from Lemma \ref{lem:quasiconformal}-(iii), and the proof of Theorem \ref{th:fun} is done.
\end{proof}

\begin{corollary} \label{co:fun}
Let $X=(X_j)_{j=1,2,3}:M \to \r^3$ be a harmonic immersion, where $M$ is an open Riemann surface, and let $\Phi=(\Phi_j)_{j=1,2,3}$ denote its Weierstrass data. Then the following statements are equivalent:
\begin{enumerate}[\rm (i)]
\item $X$ is complete and of FTC.
\item $X$ is complete and algebraic, and the Gauss map $\Ggot:M \to \s^2$ of $X$ extends continuously to the compactification $\overline{M}$ of $M.$
\item $X$ is algebraic, and for each end $P \in \overline{M}-M,$ $I_P \geq 1$ and there exists  $R_P \in O(3,\r)$ such  that $$\text{${\rm Ord}_P\left((R_P \circ \Phi)_3\right)<{\rm Ord}_P\left((R_P \circ \Phi)_2\right)={\rm Ord}_P\left((R_P \circ \Phi)_1\right)=I_P+1$ and $\frac{(R_P \circ \Phi)_2}{(R_P \circ \Phi)_1}(P) \notin \r.$}$$ 
\item $X$ is complete, algebraic and QC.
\item $X$ is complete, QC and  $\int_M |\Kcal| dS<+\infty$.   
\end{enumerate}
Furthermore, if the above statements hold then $X$ is proper, $\|\sigma\|^2$ is bounded and $X(M)$ viewed from infinity looks like a finite
collection of flat planes (with multiplicity) that pass through the origin.
\end{corollary}
\begin{proof} In either case $M$ has finite topology (take into account Huber theorem \cite{huber} in cases (i) and (v)).  The remaining follows from Theorem \ref{th:fun} applied to each annular end of $M.$
\end{proof}

In the sequel, if $X:M \to \r^3$ is a complete harmonic immersion of FTC, we will still denote by $\Ggot$ and $g$ the continuous extensions of the Gauss map and complex Gauss map of $X$ to $\overline{M}.$

The following corollary compiles some Osserman \cite{O1} and Jorge-Meeks \cite{JM} results for complete harmonic immersions of FTC.
\begin{corollary} \label{co:gaussmap}
Let $X$ be a complete non-flat harmonic immersion with FTC of an open Riemann surface $M$ into $\r^3,$ and write $M=\overline{M}-\{Q_1,\ldots,Q_k\}.$  Then the following assertions hold:
\begin{enumerate}[{\rm (a)}]
\item The Gauss map  $\Ggot:\overline{M} \to \s^2$ of $X$ is a finitely branched analytical covering. 
\item $Q_j$ is a branch point of $\Ggot$ of order $I_{Q_j}-{\rm Ord}_{Q_j}(\langle \Phi,G(Q_j)\rangle),$ where $\Phi=(\Phi_j)_{j=1,2,3}$ are the Weierstrass data of $X,$ $j=1,\ldots,k.$ 
\item If we write $\nu$  and $Deg(\Ggot)$ for the genus of $\overline{M}$ and the degree of $\Ggot,$ respectively, then
$$\int_M \Kcal dS= -4\pi  Deg(\Ggot)=-2 \pi\big(2 \nu-2+ \sum_{j=1}^k (I_{Q_j}+1)\big).$$
\end{enumerate}
\end{corollary}
\begin{proof}  $X$ is QM by Lemma \ref{lem:quasiconformal}-(ii), hence the topological part of item (a) follows from  \cite{O1}. (An alternative proof of this fact  is also implicit in the following arguments.)

 Let us see both that $\Ggot$ is analytic and item (b) (for the first part, it only remains to check that $\Ggot$ is analytic at the ends, see Lemma \ref{lem:gauss}).   By Lemma \ref{lem:gauss}, we only have to study the local behavior of $g$ around the ends of $M.$ Take $Q \in \{Q_1,\ldots,Q_k\}$ and fix an annular end $A\cong \overline{\d}-\{0\}$ corresponding to the end $Q.$ Up to a suitable linear transformation, we can suppose  that 
$$\Phi_1(z)=\frac{1+z h_1(z)}{z^{n+m}}dz,\quad \Phi_2(z)=\frac{\imath+z h_2(z)}{z^{n+m}}dz\quad\text{and}\quad \Phi_3(z)=\frac{c+z h_3(z)}{z^n}dz,$$
where $n,$ $m\in\n,$ $c\in\c-\{0\}$ and $h_j:\overline{\d} \to \c$ is holomorphic $\forall j=1,2,3.$

Therefore $g(Q)=0$ and by equation \eqref{eq:gauss} $g(z)= c z^m (1+\Ocal(z)),$ where both the real and imaginary part of $\Ocal(z)$ are real analytic functions and $\limsup_{z \to 0} |\Ocal(z)/z|<+\infty.$ In particular, $\Ggot$ is analytic at $Q$ and $g(\overline{A})$ is a neighborhood of $0.$ Moreover,   $g$ is regular at $Q$ if and only if $m=1,$ and when $m>1$ the map $g$ has an (isolated) topological branch point at $Q$ of order $m-1\geq 1.$   Since $Q$ is an arbitrary end of $M,$ item (b) holds.

To finish, let us check (c). By Osserman results \cite{O1}, $\int_M \Kcal dS= -4\pi  Deg(\Ggot).$ Up to a new rigid motion, assume that $g(Q_j) \neq 0,\infty$ for all $j=1,\ldots,k,$ and $dg (P)\neq 0$ for all $P \in g^{-1}(\{0,\infty\}).$  By Lemma \ref{lem:gauss}-(b), $\Phi_3$ has exactly $2 Deg(\Ggot)$ zeros and $\sum_{j=1}^k (I_{Q_j}+1)$ poles on $\overline{M}.$ By the Riemann relation, $2 Deg(\Ggot)-\sum_{j=1}^k (I_{Q_j}+1)=2 \nu-2$ and we are done.
\end{proof}

Osserman \cite{O2} proved that the normals to a complete non-flat minimal surface of FTC in $\r^3$ assume every spherical value with at most three exceptions (two in the genus zero case). Furthermore, he also proved that the normals to a complete QM surface of FTC can omit at most four values unless the surface is a plane (see \cite{O1}). 
 As a consequence of Corollaries \ref{co:fun} and \ref{co:gaussmap}, these results hold for complete harmonic immersions of FTC in $\r^3$ as well. The proof of the following corollary  is essentially that given by Osserman (see \cite{O1,O2,O}).
\begin{corollary} \label{co:4puntos}
Let $M$ be an open Riemann surface and let $X:M \to \r^3$ be a complete harmonic immersion of FTC. If the Gauss map of $X$ omits more than 3 spherical values, then $X(M)$ is a plane. The same conclusion holds when $M$ is a planar domain and the Gauss map of $X$ omits more than two values.
\end{corollary}

It still remains an open question to determine the precise number of values that the normals to a complete non-flat minimal surface of FTC in $\r^3$ can omit. Since the Gauss map of the minimal catenoid omits two directions then the precise bound must be ``two'' or ``three''. This problem naturally extends to complete non-flat harmonic immersions with FTC, and still remains open even in this more general setting. 

\subsection{Finite planes}

In this subsection we will prove some Picard type theorems for harmonic immersions (related results for minimal surfaces can be found in \cite{L}). To be more precise, we study how the existence of finite planes for a harmonic immersion influences its global geometry. Our techniques  are inspired by the arguments of Xavier in \cite{xavier3} and Meeks-P\'{e}rez results in \cite{MP}.

Let $X:M \to \r^3$ be a harmonic immersion. A plane $\Pi\subset \r^3$ is said to be {\em finite} for $X$ if $X^{-1}(\Pi)$ has only a finite number of connected components, and the surfaces $\Pi$ and $X(M)$ meet transversally at all but at most finitely many points of  $X(M) \cap \Pi.$ 

\begin{figure}[ht]
    \begin{center}
    \scalebox{0.35}{\includegraphics{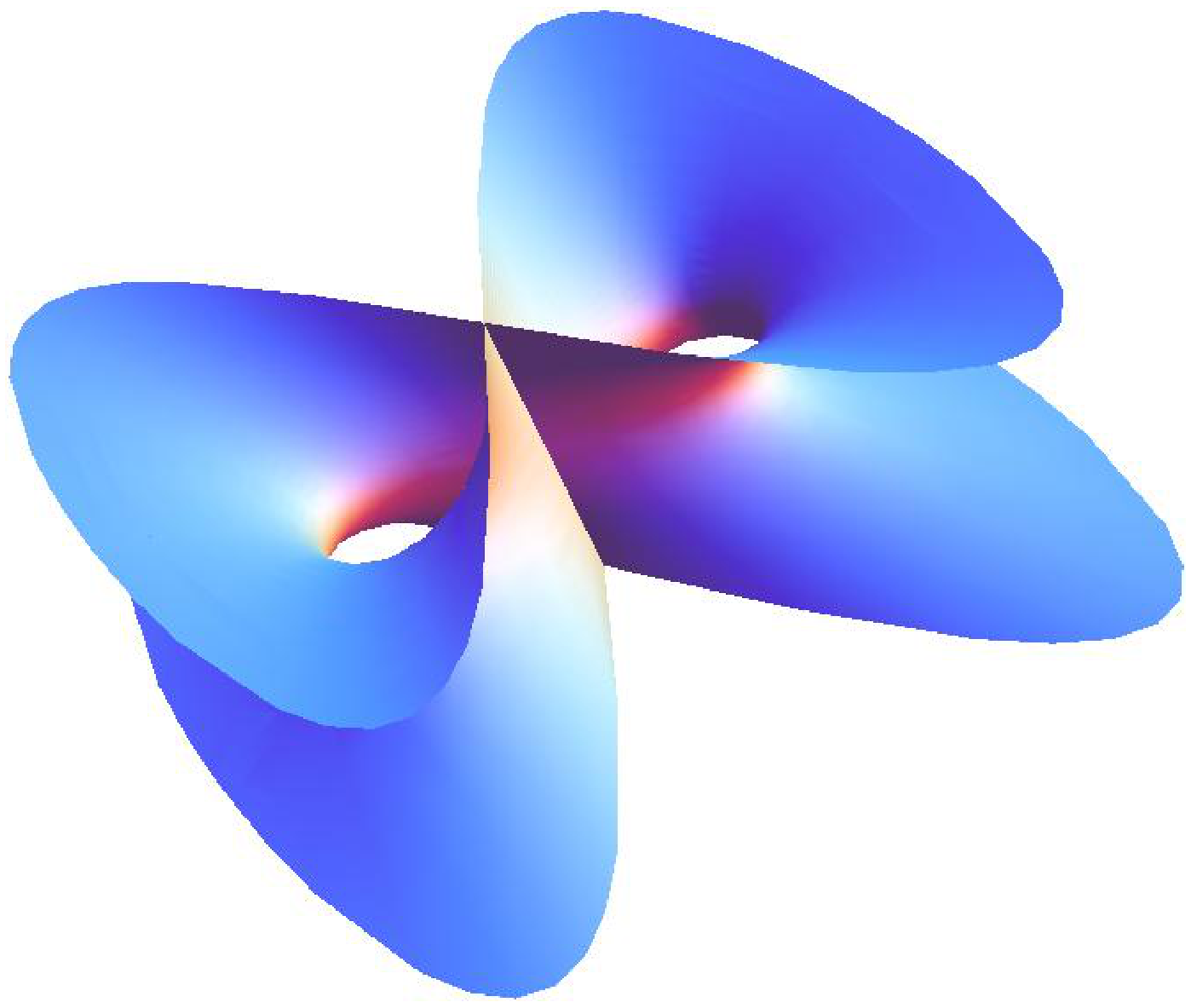}}
        \end{center}
\caption{The surface $Y(\c-\{0\}).$}\label{fig:rara}
\end{figure}

\begin{theorem}\label{th:anular}
Let $A=\{r<|z|\leq 1\},$ $r\in [0,1[,$  be an annular end, and let $X=(X_j)_{j=1,2,3}:{A} \to \r^3$ be a complete harmonic  immersion. The following statements hold:
\begin{enumerate}[\rm (a)]
\item $X$ is of FTC if and only if  $X$ is QC, $\sup_{A} \|\sigma\|^2<+\infty$ and $X$ has two non parallel finite planes.
\item $X$ is of FTC if and only if  $X$ is  QC,  proper and has two non parallel finite planes.
\item If $\sup_{A} \|\sigma\|^2<+\infty$ and $X$ has three finite planes in general position then $X$ is algebraic.
\item If $X$ is proper and has three finite planes in general position then $X$ is algebraic.
\end{enumerate}
\end{theorem}
\begin{proof} 
Let us show (a) and (b).

Assume that $X$ is QC and has two non parallel finite planes $\Pi_1$ and $\Pi_2.$ Up to a rigid motion suppose that $\Pi_1 \cap \Pi_2$ is a vertical straight line. 
\begin{claim}
If either $\sup_{A} \|\sigma\|^2<+\infty$ or $X$ is proper  then $r=0.$
\end{claim}
\begin{proof} Reason by contradiction an assume that $r>0.$ 
Put $\Pi_j=\{a_j x_1+b_j x_2=c_j\}$ and  $Y_j=a_j X_1+b_j X_2,$ $j=1,2.$ Label $\Omega = {A}-\{\Im(z)=0,\; \Re(z)>0\}$ and write $F_j:\Omega \to \c$ for the holomorphic function $Y_j +\imath Y_j^*,$ where $Y_j^*$ is the harmonic conjugate of $Y_j,$ $j=1,2.$   Xavier results \cite{xavier3} imply that $F_j$ has non-tangential limits almost everywhere. This means that for almost every $\xi \in \partial (\Omega)$ and any triangle $T_\xi\subset \Omega$ with vertex $\xi,$  the limit $\lim_{z \to \xi,\; z \in T_\xi} F_j(z)$ exists in $\overline{\c}.$  Moreover,  by the Privalov uniqueness theorem  the non-tangential limits of $F_j$ are finite  for almost every $\xi\in \partial(\Omega),$  $j=1,2.$ Combined results by Marcinkiewicz, Zygmund and Spencer, see \cite{Ko, xavier1}, imply that  $\int_{T_\xi} |\partial_z Y_j|^2 |dz|^2<+\infty,$ $j=1,2,$ and so
\begin{equation}\label{eq:areafin}
\int_{T_\xi} |\partial_z X_j|^2 |dz|^2<+\infty
\end{equation}
for almost every $\xi$ and any $T_\xi$ as above, $j=1,2.$

Label $C=\sup_{A} \frac{|\Hgot|}{\|\Phi\|^2},$ and notice that $C<1$ by quasiconformality, see Proposition \ref{pro:index}. Then 
\begin{equation} \label{eq:areafin3}
|\partial_z X_3|^2\leq \frac{C+1}{1-C} \left(|\partial_z X_1|^2 +|\partial_z X_2|^2 \right),
\end{equation}
hence from equation \eqref{eq:areafin},  
$\int_{T_\xi} |\partial_z X_3|^2 |dz|^2<+\infty$  for almost every $\xi$ and $T_\xi$ as above. Again by Marcinkiewicz, Zygmund, Spencer and Privalov  results one has that $F_3:\Omega \to \c$ has finite non-tangential limits for almost every $\xi\in \partial(\Omega),$  where $F_3=X_3+\imath X_3^*.$  

Therefore, there is   $\xi \in \partial(\Omega) \cap \{|z|=r\}$ and $T_\xi$ such that the limit  $\lim_{z \to \xi,\; z \in T_\xi} (F_j(z))_{j=1,2,3}$  exists. This is absurd when $X$ is proper. In case  $\sup_{A_0} \|\sigma\|^2<+\infty,$ a contradiction can then be reached by considering the
point $p=\lim_{z \to \xi,\; z \in T_\xi} X(z)$ in $\r^3$ and using the geometric lemma from \cite{xavier2,xavier1}. This proves the claim.
\end{proof}

Assume in addition that either $\sup_{A} \|\sigma\|^2<+\infty$ or $X$ is proper. Then, the claim gives that ${A}=\overline{\d}-\{0\}.$ By Meeks-P\'{e}rez results \cite{MP}, the function $\partial_z Y_j$ has a meromorphic extension to $\overline{\d},$ $j=1,2,$ hence the same holds for $\partial_z X_j,$ $j=1,2.$ By equation \eqref{eq:areafin3}, $\partial_z X_3$  has a meromorphic extension to $\overline{\d}$ as well,  and Theorem \ref{th:fun} shows that $X$ has FTC. 

Obviously, all the planes in $\r^3$ are finite for any algebraic harmonic immersion, and so, for any complete harmonic immersion with FTC.
The converse in (a) and (b) follow from  Theorem \ref{th:fun}.

The proof of items (c) and (d) is similar.  
\end{proof} 
The rotational harmonic horn and the example in Remark \ref{re:contra} are counterexamples for the converses in  (c) and (d), respectively. Moreover, the hypothesis of being QC plays a role in (a) and (b). Indeed, the map 
$$Y:\c-\{0\} \to \r^3,\quad Y(z)=\big( \Im(\frac{-1 + z^4}{2 z^2}), \Re(1/z + z), \log |z| \big)$$ is a non-QC proper algebraic harmonic immersion with bounded $\|\sigma\|^2$ (and of course, with infinite total curvature by Corollary \ref{co:fun}). See Figure \ref{fig:rara}.

\section{Harmonic embeddings of finite total curvature}\label{sec:embe}
In this section we  study  the global geometry of complete harmonic embeddings with FTC. We will  describe some examples of harmonic embeddings with genuine geometrical properties in contrast  to the case of minimal surfaces in $\r^3.$

We start with the following Schoen type expansion  for embedded harmonic ends of FTC (for the conformal case see \cite[Proposition 1]{schoen}).
\begin{lemma} \label{lem:embe}
Let $A$ be a conformal annulus biholomorphic to $\overline{\d}-\{0\}$, and let  $X=(X_j)_{j=1,2,3}:A \to \r^3$ be a complete harmonic immersion of FTC. 

Then,  $X$ is an embedding outside a compact set if and only if $I_Q=1$ (see Definition \ref{def:iq}), where $Q\in \overline{A}$ is the topological end of $A.$ In this case, up to removing a compact subset of $A$  and up to a linear transformation, one has $$X(A)=\{(x,u(x))\in \r^2\times \r\,|\; \|x\| \geq C\},$$ where $u:\{x \in \r^2\,|\;\|x\| \geq C\} \to \r$ is an analytic function satisfying
\begin{equation} \label{eq:expansion}
 \text{either}\;\lim_{x \to \infty} u(x)=+\infty \; \text{and}\; e^{u(x)}=\|x-(a_1,a_2) u(x)\|+\Ocal(1),\quad \text{or} \quad \|x\| u(x)=\Ocal(1).
\end{equation}
The expression $\Ocal(1)$ is used to indicate a term which is bounded on $\{\|x\|\geq C\}.$
\end{lemma}
\begin{proof} By Theorem \ref{th:fun}, $X$ is an embedding outside a compact subset of $A$ if and only if $I_Q=1.$ 

Let $\Phi=(\Phi_j)_{j=1,2,3}$ denote the Weierstrass data of $X.$ By Theorem \ref{th:fun} and up to a rigid motion, there exists a conformal parameter $z$ on $A$ so that $z(Q)=0$ and  
$$\Phi_1=\big(\frac{1}{z^2}-\frac{a_1}{z}+h_1(z)\big)dz,\;\; \Phi_2=\big(\frac{\imath}{z^2}-\frac{a_2}{z} +h_2(z)\big) dz,\;\;\Phi_3=\big(-\frac{\delta}{z}+h_3(z)\big)dz,$$ where   $a_1,$ $a_2 \in \r,$ $\delta\in \{0,1\}$ and $h_j$ is bounded, $j=1,2,3.$ So $X_3=-\delta \log|z|+\Ocal(|z|),$ 
$$X_1-a_1X_3=-\frac{\Re(z)}{|z|^2}+(\delta-1) a_1 \log |z| +\Ocal(|z|),\;\; X_2-a_2X_3=\frac{\Im(z)}{|z|^2}+(\delta-1) a_2 \log |z| +\Ocal(|z|),$$ 
where $\Ocal(|z|)$ is used to indicate a term which is bounded by a constant times $|z|.$
The asymptotic expansions \eqref{eq:expansion} easily follow from these expressions.
\end{proof}
An embedded annular end $X=(X_j)_{j=1,2,3}:A \to \r^3$ of FTC is said to be {\em normalized}  if $X(A)=\{(x,u(x))\,|\; \|x\| \geq C\}$ and \eqref{eq:expansion} holds. In this case the limit tangent plane at the end is horizontal.

Let $Y:A\to\r^3$ be an embedded annular end of FTC, write $A=\overline{A}-\{Q\},$ and take a linear transformation $R \in Gl(3,\r)$ so that $X=(X_j)_{j=1,2,3}:=R \circ Y$ is normalized.

 $Y$ is said to be a {\em harmonic planar end} if  $\|(X_1,X_2)\| X_3=\Ocal(1).$  In the proof of Corollary \ref{co:gaussmap} it is shown that $Y$ is planar if and only if $Q$ is a branch point of  the Gauss map  $\Ggot$ of $Y.$ When $Y$ is planar and the branching number of $\Ggot$  at $Q$ is one,  then $X$ is said to be a harmonic  planar end  of {\em Riemann type}. In this case $\Phi_ 3=\partial_z X_3 (Q)\neq 0,\infty,$ hence it is easy to see that, up to removing a compact subset of $A,$
\begin{equation} \label{eq:embe}
\text{the limit tangent plane of $Y$ at $Q$ divides $Y(A)$ into exactly two connected components.} 
\end{equation}

$Y$ is said to be a {\em harmonic catenoidal end} if $e^{X_3}=\|(X_1,X_2)-(a_1,a_2) X_3\|+\Ocal(1).$ In this case, the  straight line  $R^{-1}\left( \{(a_1,a_2,1)t\,|\; t \in \r\}\right)$ is said to be the {\em axis} of $Y(A).$ Note that the axis and the normal vector at infinity of a harmonic catenoidal end are not necessarily parallel. This only happens when the associated normalized end is asymptotic to a rotational harmonic  half catenoid ($a_1=a_2=0$ in Lemma \ref{lem:embe}), or equivalently, if $X$ has vertical flux.

The following corollary follows from Corollary \ref{co:fun}, Lemma \ref{lem:embe} and well known arguments by Jorge-Meeks \cite{JM}.

\begin{corollary} \label{co:ordering}
Let $X=(X_j)_{j=1,2,3}$ be a complete harmonic immersion with FTC of an open Riemann surface $M$ into $\r^3.$ Then the following statements hold:
\begin{enumerate}[{\rm (a)}]
\item If $M$ has an only end and $X$ is an embedding outside a compact subset, then $X(M)$ is a plane.
\item If $X$ is an embedding and $M$ has more than one end, then all its ends are either catenoidal or planar, and up to a rigid motion, the limit tangent planes of the ends are parallel to  $\{x_3=0\}.$ Furthermore, ordering the ends by their heights over $\{x_3=0\},$ two consecutive ends have opposite normal, and the top and bottom ends are of catenoidal type and with logarithmic growth of opposite sign. 
\end{enumerate} 
\end{corollary}
\begin{proof}
In case (a), the equation \eqref{eq:expansion} shows that, up to a rigid motion,  $X_3$ is positive, hence constant by the maximum principle (recall that $M$ is parabolic, see Corollary \ref{co:fun}). 

To check (b), use Lemma \ref{lem:embe} to infer that all the ends of $X$ are parallel. Note that the properly embedded surface $X(M)$ splits $\r^3$ into two connected components, and  the Gauss map of $X$ points to one of them. Since $M$   carries no non-constant positive harmonic functions,  the top and bottom ends of $X$ must be of catenoidal type and with logarithmic growth of opposite sign. 
\end{proof}

\subsection{The harmonic catenoids}
In this subsection we classify all complete harmonic embeddings with FTC of an open annulus in $\r^3.$

Consider the Weierstrass data 
\begin{equation} \label{eq:wei-cat}
M=\c-\{0\},\;\; \Phi_1=\left(\frac{1}{z^2}+\frac{r_1}{z}+\frac{\alpha+\overline{\beta}}{2} \right) dz,\;\; \Phi_2=\left(\frac{\imath}{z^2}+\frac{r_2}{z}+\frac{\alpha-\overline{\beta}}{2\imath} \right) dz,\;\; \Phi_3=\frac{dz}{z},
\end{equation}
where $\alpha,$ $\beta \in \c,$ $|\alpha|\neq |\beta|$ and $r_1,$ $r_2 \in \r.$ Set 
$$X_{\alpha,\beta}:\c-\{0\} \to \r^3, \quad X=\Re \int \left(\left(\frac{1}{z^2}+\frac{r_1}{z}+\frac{\alpha+\overline{\beta}}{2}  \right) dz,\left(\frac{\imath}{z^2}+\frac{r_2}{z}+\frac{\alpha-\overline{\beta}}{2 \imath}  \right) dz,\frac{dz}{z} \right) .$$

Call $\Omega=\c^2-\big(\{(u,v) \in \c^2\,|\; |u|\leq |v|\} \cup \{(u,v) \in \c^2\,|\;\Re (u)\geq 0\;\text{and}\; |\Im (u)|\leq |v|\}\big).$ 

\begin{lemma} \label{lem:cat-immer}
The harmonic map $X_{\alpha,\beta}$ is well defined, and it is an immersion if and only if $(\alpha,\beta) \in \Omega.$
\end{lemma}
\begin{proof} Since $r_1,$ $r_2 \in \r,$ $\Phi$ has no real periods and $X$ is well defined.

By Lemma \ref{lem:wei}, $X$ is an immersion if and only if $\frac{|\sum_{j=1}^3 \Phi_j^2|}{\sum_{j=1}^3 |\Phi_j|^2}<1$ on $\c-\{0\}.$ Since $dz/z$ is holomorphic and never vanishes on $\c-\{0\},$ this inequality is equivalent to the fact that
$$(\Phi_1/\Phi_3,\Phi_2/\Phi_3)\notin \r^2\;\; \text{everywhere on $\c-\{0\}$}.$$

Assume for a moment that $(\Phi_1/\Phi_3,\Phi_2/\Phi_3)(w)\in \r^2$ for some $w \in \c-\{0\}.$ Therefore, 
$$\frac{\alpha+\overline{\beta}}{2}  w-\frac{\overline{\alpha}+\beta}{2}  \overline{w}+\frac{1}{w}-\frac{1}{\overline{w}}=0\;\;\text{and}\;\;
\frac{\alpha-\overline{\beta}}{2\imath}  w+ \frac{\overline{\alpha}-\beta}{2\imath} \overline{w}+\frac{\imath}{w}+\frac{\imath}{\overline{w}}=0,$$ hence
\begin{equation} \label{eq:c-i-c}
\overline{w} \left(\alpha w-\beta \overline{w} \right)=2.
\end{equation}
If we set $f:(0,\infty) \times [0,2\pi] \to \c,$ $f(t,\theta)=t^2 \left(\alpha-\beta e^{-2 i \theta} \right),$ equation \eqref{eq:c-i-c} says that
$X$ is an immersion if and only if
$2 \notin f((0,\infty) \times [0,2\pi]),$ that is to say, if and only if $2/t^2\notin C_{\alpha}(|\beta|)$ for all $t>0,$  where $C_{\alpha}(|\beta|)\subset \c$ is the circle of center $\alpha$ and radius $|\beta|.$

On the other hand,  $2/t^2\in C_{\alpha}(|\beta|)$ for some $t>0$ if and only if either $|\alpha|<|\beta|$ or $\Re(\alpha)>0$ and $|\Im(\alpha)|\leq |\beta|.$ Taking into account that $|\alpha|\neq|\beta|$  we are done. 
\end{proof}

\begin{definition}
For any $(\alpha,\beta)\in\Omega,$ the harmonic immersion $X_{\alpha,\beta}:\c-\{0\}\to\r^3$ is said to be a {\em harmonic catenoid}.
\end{definition}

\begin{theorem} \label{th:cat-todas}
For any $(\alpha,\beta) \in \Omega,$ $X_{\alpha,\beta}$  is a complete harmonic embedding of total curvature $-4 \pi.$

Furthermore, if $X$ is a complete harmonic embedding with  FTC of an open conformal annulus $M$ into $\r^3,$ then up to biholomorphisms and linear transformations $X=X_{\alpha,\beta}$ for some $(\alpha,\beta)\in \Omega.$
\end{theorem}

\begin{proof} Let us check that $X_{\alpha,\beta}$ is a complete embedding of FTC for all $(\alpha,\beta) \in \Omega.$ By Lemma \ref{lem:cat-immer} and Corollaries \ref{co:fun} and \ref{co:gaussmap}, $X_{\alpha,\beta}$ is a complete harmonic immersion of total curvature $-4 \pi.$

Integrating \eqref{eq:wei-cat} and using polar coordinates $z=m e^{\imath t}$, one has
\begin{multline*}
X_{\alpha,\beta}(m,t)=\left(\frac{(-2+\Re(\alpha+\beta) m^2)\cos(t)+\Im(\beta-\alpha) m^2 \sin(t)+ 2r_1 m \log(m)}{2 m},\right. \\
\left. \frac{(-2+\Re(\alpha-\beta)m^2) \sin(t)+ \Im(\alpha+\beta) m^2 \cos(t)+2r_2 m \log(m)}{2 m},\log(m) \right).
\end{multline*}
For the sake of simplicity, write $X_{\alpha,\beta}=(X_j)_{j=1,2,3}.$ Assume that $X_{\alpha,\beta}(m_1,t_1)=X_{\alpha,\beta}(m_2,t_2).$ Note that $m_2=m_1$ and look the   
 indentity $(X_1,X_2)(m_1,t_1)=(X_1,X_2)(m_1,t_2)$ as a system of linear equations with variables $\cos(t_1)$ and $\sin(t_1).$ The determinant of the matrix of this system is 
 $-\Re(\alpha)+\frac{1}{m^2}+\frac{m^2}{4} (|\alpha|^2-|\beta|^2),$ which is positive for $(\alpha,\beta) \in \Omega.$ Therefore, the only solution is $(\cos(t_1),\sin(t_1))=(\cos(t_2),\sin(t_2)),$ and so $e^{\imath t_1}=e^{\imath t_2}.$ This proves that $X_{\alpha,\beta}$ is injective, hence an embedding (recall that $X_{\alpha,\beta}$ is proper).

For the second part of the theorem, consider a conformal annulus $M$ and   a complete harmonic embedding $X=(X_j)_{j=1,2,3}:M \to \r^3$ of FTC. Write $(\Phi_j)_{j=1,2,3}=\partial_z X$ for the Weierstrass data of $X.$  Corollary \ref{co:ordering} implies that $X$ has two catenoidal ends with opposite logarithmic growths. Furthermore, up to biholomorphisms and linear isometries,  Theorem \ref{th:fun} (see also Lemma \ref{lem:embe}) gives that 
$$M=\c-\{0\},\;\; \Phi_1=\left(\frac{1}{z^2}+\frac{r_1}{z}+a_1 \right) dz,\;\; \Phi_2=\left(\frac{\imath}{z^2}+\frac{r_2}{z}+a_2 \right) dz,\;\; \Phi_3=\frac{1}{z},$$
where $a_1,$ $a_2 \in \c-\{0\},$ $\Im(a_1 \overline{a}_2)\neq 0$ and $r_1,$ $r_2 \in \r.$ For the condition on $r_1$ and $r_2,$ take into account that $X$ has no real periods. Write $\alpha=a_1+\imath a_2$ and $\beta= \overline{a}_1+\imath \overline{a}_2,$ and notice that  $\Im(a_1 \overline{a}_2)\neq 0$ if and only if $|\alpha| \neq |\beta|.$ By Lemma \ref{lem:cat-immer}, $X \in  \{X_{\alpha,\beta}\,|\;(\alpha,\beta) \in \Omega\}$ and we are done.
\end{proof}

\begin{remark} $X_{\alpha,\beta}$ has vertical flux if and only if $r_1=0=r_2.$ A straightforward computation gives that
$X_{\alpha,\beta}$ is linearly equivalent to a rotational harmonic catenoid if and only if $X_{\alpha,\beta}$ has vertical flux and $\beta=0.$ Furthermore, in this case  $X_{\alpha,\beta}$ is in fact rotational. 

In Figure \ref{fig:cat-flux} can be found the pictures of $X_{\alpha,\beta}(\c-\{0\})$ for $(r_1,r_2,\alpha,\beta)=(2,0,-3+3 \imath,-1-\imath)$ on the left, and for  $(r_1,r_2,\alpha,\beta)=(0,0,-3+3 \imath,-1-\imath)$ on the right.
\end{remark}

\subsection{Genus zero harmonic embeddings with vertical flux}
It is well known that there are no properly embedded minimal surfaces in $\r^3$ of genus zero, finite topology and more than two ends \cite{lop-ros,collin}. Furthermore, the minimal catenoid is the only complete non-flat embedded  minimal surface of FTC with vertical flux \cite{lop-ros,perez-ros}. Contrary to the minimal case, in this subsection we construct harmonically embedded triply-connected planar domains in $\r^3$ with FTC and vertical flux.

Set $M=\c-\{-1,1\},$ and for each $b$ and $c \in \r$ consider the harmonic map
$$X:M \to \r^3, \quad X(z)=\left(\Im\big(z (b-2+\frac8{z^2-1}) \big),\Re\big(z (c-2+\frac8{z^2-1}) \big), 6\log \big|\frac{z-1}{z+1}\big| \right).$$
If we write $(\Phi_j)_{j=1,2,3}=\partial_z X,$ a straightforward computation gives that:
$$\Phi_1=-\imath \frac{6+b +(12-2 b) z^2+(b-2) z^4}{(z^2-1)^2}dz,\;\;  \Phi_2=\frac{6+c +(12-2 c) z^2+(c-2) z^4}{(z^2-1)^2}dz,\;\; \Phi_3=\frac{12dz}{z^2-1}.$$

\begin{lemma} \label{lem:flujo}
If  $b>3$ and $c<2$ then $X$ is a harmonic immersion.
\end{lemma}
\begin{proof} Notice that $\Phi_3$ never vanishes on $M.$ As in the preceding theorem, it suffices to prove that $(f_1,f_2) \notin \r^2$ everywhere on $M,$ where  $f_1=\frac{\Phi_1}{\Phi_3}$ and $f_2=\frac{\Phi_2}{\Phi_3}.$ For the sake of simplicity, write $t=z^2$ and $h_j(t)=f_j(\sqrt{t}),$ $j=1,2.$ Observe that
$$h_1(t)= -\imath \frac{6+b +(12-2 b) t+(b-2) t^2}{12(t-1)},\quad h_2(t)= \frac{6+c +(12-2 c) t+(c-2) t^2}{12(t-1)}.$$

Let us see that $(h_1,h_2)(t) \notin \r^2$ when $t \in \r.$ Indeed, in this case $h_1(t) \in \imath \r,$ hence $h_1(t) \in \r$ if and only if $h_1(t)=0.$ This is not possible since $b>3.$ 

Reason by contradiction and assume there exists $t \in M-\r$ such that $(h_1,h_2)(t) \in \r^2.$ As $h_2(t) \in \r,$ then an elementary computation gives that 
$\overline{t}=\frac{18 - c - 2 t + c t}{(-2 + c) (-1 + t)}.$ This equality and the fact that $h_1(t) \in \r$ imply that 
$$-24 + 14 b - 6 c + b c + c^2 -  2 (-2 + c) (-12 + b + c) t + (-2 + c) (-4 + b + c) t^2=0.$$ As $t$ is not real, the discriminant $(2 - c) (24 + b^2 + 2 b (-4 + c) + (-12 + c) c)$ must be negative, which contradicts that $c<2$ and $b>3.$
\end{proof}

\begin{theorem} \label{th:flujo}
For any $b>3$ and $c<2,$ $X$ is a complete harmonic embedding with FTC and vertical flux. 
\end{theorem}
\begin{proof}By Lemma \ref{lem:flujo} and Corollaries \ref{co:fun} and \ref{co:gaussmap}, $X$ is complete, of FTC and the degree of its Gauss map $\Ggot:\overline{\c}\to\s^2$ is $Deg(\Ggot)=2.$ On the other hand, clearly $X_j$ is the real part of a holomorphic function $\Xcal_j:M\to\c,$ $j=1,2,$ hence $X$ has vertical flux.

Let us check that $X$ is an embedding. Observe that $X(-z)=-X(z),$ $X(\overline{z})=(-X_1,X_2,X_3)(z),$ $X(\imath \r)\subset\{x_2=x_3=0\}$ and $X(\r-\{-1,1\})\subset \{x_1=0\}.$ 
Therefore, $\Ggot(\r \cup \imath \r \cup \{\infty\}) \subset \s^2 \cap \{x_1=0\}.$ Since $\Ggot$ has degree two, then $\Ggot^{-1}(\s^2 \cap \{x_1=0\})=\r \cup \imath \r \cup \{\infty\},$ and $0$ and $\infty$ are the only ramification points of $\Ggot.$ Consequently  $\Ggot(M_{i,j})\subset \s^2\cap \{(-1)^{i+j} x_1 > 0\},$ where  $M_{i,j}$ is the quadrant $\{z \in M\,|\; (-1)^{i}\Re(z) > 0,\; (-1)^{j}\Im(z) > 0\},$  $i,j \in \{0,1\}.$ 

Taking into account the geometry of harmonic catenoidal ends and Riemann type harmonic planar ends (see Lemma \ref{lem:embe} and  \eqref{eq:embe}), we infer that $X(M_{0,j})-\{x_2=x_3=0\}$ is a graph over the planar domain $\Omega_0$ in $\{x_1=0\}$ bounded by the Jordan arcs $X([0,1)),$ $X((1,+\infty))$ and the half line $\{x_3=x_1=0,\; x_2 \leq 0\}$ (see Figure \ref{fig:omega}) and 
$X|_{\overline{M}_{0,j}}$ is injective, $j=1,2.$ 
 Likewise, $X(M_{1,j})-\{x_2=x_3=0\}$ is a graph over the planar domain $\Omega_1:=-\Omega_0$ in $\{x_1=0\},$ and $X|_{\overline{M}_{1,j}}$ is injective, $j=1,2.$ Notice that $\Omega_1$ is  bounded by  $X((-1,0]),$ $X((-\infty,-1))$ and the half line $\{x_3=x_1=0,\; x_2 \geq 0\},$ and $\Omega_0\cap \Omega_1=\emptyset.$ 
\begin{figure}[ht]
    \begin{center}
    \scalebox{0.75}{\includegraphics{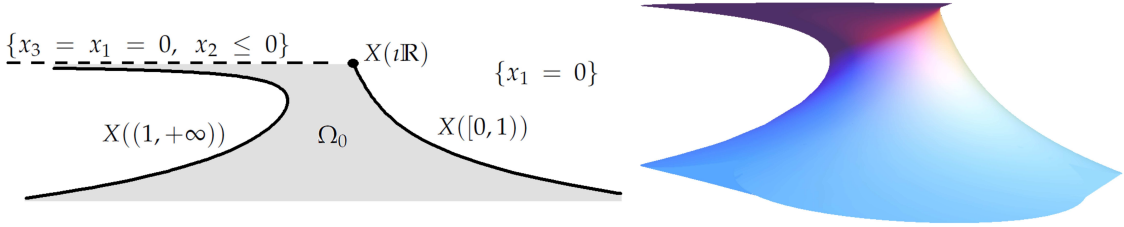}}
        \end{center}
\caption{The planar domain $\Omega_0\subset\{x_1=0\}$ and $X(M_{0,j})-\{x_2=x_3=0\}$ as a graph over it.}\label{fig:omega}
\end{figure}

As $(-1)^j X_1>0$ and $(-1)^i X_3 <0$  on $M_{i,j}$ for all $i,j$ and $M=\cup_{i,j} \overline{M}_{i,j},$ then $X$ is injective. This concludes the proof.
\end{proof}

The picture on the left in Figure \ref{fig:flujo+toro} corresponds to $X(\c-\{-1,1\})$ for $b=4$ and $c=0.$

\subsection{The harmonic catenoidal tori} 
Numerical evidences of the existence of harmonically embedded tori with two catenoidal ends were pointed out by Weber in \cite{we} (see the picture on the right in Figure \ref{fig:flujo+toro}). There are no equivalent surfaces in the minimal case. In fact, the catenoid is the only complete embedded  minimal surface of finite genus with two topological ends  \cite{schoen,collin,colding-mini}. Here we present a rigorous proof of the regularity and embeddedness of these surfaces.

For each $a \in (0,1)$ and $b \in \r,$ set  
$$\overline{M}_a=\{(z,w)\in \overline{\c}^2\;|\; w^2=\frac{(z-a) (a z-1)}{z}\},\quad M_a=\overline{M}_a-\{(0,\infty),(\infty,\infty)\},$$
\begin{equation} \label{eq:wei-k}
\Phi_1=\frac{\imath (z^2-1)}{z^2 w} dz,\quad \Phi_2=\frac{z^2+b z+1}{z^2 w} dz,\quad \Phi_3=\frac{dz}{z}.
\end{equation}

We will need the conformal transformations in $\overline{M}_a$ given by:
$$J(z,w)=(z,-w),\quad T(z,w)=(1/\overline{z},-\overline{w}),\quad S(z,w)=(\overline{z},\overline{w}).$$ Observe that $J(M_a)=M_a,$ and likewise for $T$ and $S.$ 

Let us study the period problem associated to the above Weierstrass data. First, notice that $\Phi_1= \frac{2\imath}{a} d w$ is exact on $M_a$ and $\Phi_3$ has no real periods on $M_a.$ As Weber showed in \cite{we}, for any $a \in (0,1)$ there exists a unique $b(a) \in \r$ such that the 1-form $\Phi_2$  for $b=b(a)$ has no real periods on $M_a.$ In order to prove the regularity of the arising harmonic immersion, we need a little more information.

\begin{lemma} \label{lem:b}
 $b(a)\in (-2,0).$
 \end{lemma}
\begin{proof} Consider two closed curves $c_1$ and $c_2$ in the $z$-plane illustrated in Figure \ref{fig:c}, and let $\gamma_j$ be any lifting of $c_j$ to $M_a,$ $j=1,2.$ Notice that $\gamma_j$ is a closed curve in $M_a$ as well, $j=1,2,$ and  $\{\gamma_1,\gamma_2\}$ is a basis of the homology group $\Hcal_1(\overline{M}_a,\z).$ Observe that $\Phi_2$ has no residues at neither $(0,\infty)$ nor $(\infty,\infty).$ Moreover,  since $S^*(\Phi_2)=\overline{\Phi}_2$ and $S_*(\gamma_2)=-\gamma_2$ then $\Re \big(\int_{\gamma_2} \Phi_2 \big)=0.$ 
\begin{figure}[ht]
    \begin{center}
    \scalebox{0.4}{\includegraphics{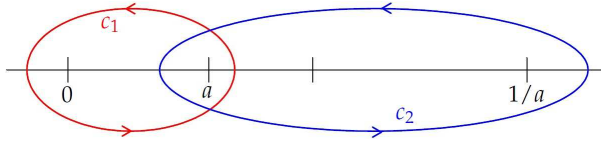}}
        \end{center}
\caption{The curves $c_1$ and $c_2$ in the $z$-plane.}\label{fig:c}
\end{figure}

To finish, it suffices to prove that the solution $b(a)$ of the equation $\int_{\gamma_1} \Phi_2=\int_{\gamma_1} \frac{z^2+1}{z^2 w} dz +b\int_{\gamma_1}\frac{1}{z w} dz=0$ lies in $(0,2).$ Indeed, recall that $\frac{z^2-1}{z^2 w}dz$ is exact. Then,  integrating by parts, $\int_{\gamma_1} \Phi_2=\int_{\gamma_1} \frac{2}{w} dz +b\int_{\gamma_1}\frac{1}{z w} dz.$ 

Notice that $\int_{\gamma_1}\frac{1}{z w} dz=2 \int_{[0,a]} \frac{1}{z w} dz \neq 0,$ where up to changing the branch  $w(z) \in \imath \r_+$ for all $z \in (0,a).$ Since    $\int_{[0,a]} \frac{1}{|w|} dz< \int_{[0,a]} \frac{1}{z |w|} dz,$  one has that $$b(a):=-2 (\int_0^a \frac{1}{z |w|} dz)/(\int_0^a \frac{1}{|w|} dz)$$ is the unique solution of $\int_{\gamma_1} \frac{2}{w} dz +b\int_{\gamma_1}\frac{1}{z w} dz=0$ and $b(a)\in (-2,0).$ 
\end{proof}
From now on $b=b(a).$
\begin{theorem} \label{th:i}
For each $a \in (0,1),$ the map $X:M_a \to \r^3,$ $X(P)=\Re(\int^P \Phi)$ is a complete harmonic AC embedding of FTC. 
\end{theorem}
\begin{proof} By Lemma \ref{lem:b}, $X$ is well defined. Let us show that $X$ is an immersion. By Lemma \ref{lem:wei}, we only have to check that $F:=\frac{|\sum_{j=1}^3 \Phi_j^2|}{\sum_{j=1}^3 |\Phi_j|^2}<1$ on $M_a.$ As $dz/w$ is holomorphic and never vanishes on $M_a,$ then the inequality $F<1$ is equivalent to the fact that 
$\{\Re(w \frac{\Phi}{dz} ),\Im(w \frac{\Phi}{dz})\}$ are $\r$-linearly independent everywhere on $M_a.$

A straightforward computation gives that this condition holds at the points $(a,0)$ and $(1/a,0),$ just the zeros of $\Phi_3$ in $M_a.$   Therefore, $X$ is an immersion if and only if  
\begin{equation} \label{eq:l-d}
\text{$(\frac{\Phi_1}{\Phi_3},\frac{\Phi_2}{\Phi_3}) \notin \r^2$ everywhere on $M_a-\{(a,0),(1/a,0)\}.$}
\end{equation} 
For the sake of simplicity, write $f_1=\frac{\Phi_1}{\Phi_3}=\frac{\imath (z^2-1)}{z w}$ and $f_2=\frac{\Phi_2}{\Phi_3}=\frac{(z^2+b z+1)}{z w}.$  

Let us see first that \eqref{eq:l-d} holds  when $|z|=1.$ Indeed, in this case $w^2=-|z-a|^2$ then $w \in \imath \r-\{0\}.$ Moreover,  $(\frac{\imath (z^2-1)}{z},\frac{(z^2+b z+1)}{z})\in \r^2-\{(0,0)\}$ (recall that $b \in (-2,0)$), and so $f_1$ and $f_2$ can not be both real numbers. 

Reason by contradiction and assume there exists $P\in M_a-\{(a,0),(1/a,0)\}$ such that $(f_1,f_2)(P)\in\r^2.$ Note first that $f_1$ and $f_2$ have no common zeros. Take $j \in \{1,2\}$ so that $f_j(P)\neq 0$ and write $\{i,j\}=\{1,2\}.$ Since $\Im(\frac{f_i}{f_j}(P))=0$ and $|z(P)|\neq 1,$ a straightforward computation gives that $|z(P)|^2+b\Re(z(P))+1=0.$  However, $|z|^2+b\Re(z)+1 >0$ for all $z \in \c$ (recall that $b\in(-2,0)$). 

This contradiction shows that $X$ is an immersion.

It is obvious that $i^X:=\limsup_{P \to \infty} F(P)=0,$ hence $X$ is AC.

By Corollaries \ref{co:fun} and \ref{co:gaussmap}, $X$ is complete and of FTC, and the degree of its Gauss map $\Ggot:\overline{M}_a\to\s^2$ is $2.$

Finally, let us prove that $X$ is an embedding. First note that the  $B=|z|^{-1}(\{1\})$ is the set of fixed point of the antiholomorphic involution $T:\overline{M}_a \to \overline{M}_a.$ Since $T^*(\Phi)=(\overline{\Phi}_1,\overline{\Phi}_2,-\overline{\Phi}_3),$ then up to a translation $X \circ T=\Tcal \circ X,$ where $\Tcal:\r^3 \to \r^3$ is the reflection about the plane $x_3=0.$ As a consequence, $X(B)\subset \{x_3=0\}$ and $\Ggot(B) \subset  \s^2 \cap \{x_3=0\}.$ A straightforward computation (see equation \eqref{eq:gauss}) gives that $\Ggot$ has no ramification points on $B.$ Since  $B$ consists of two Jordan curves $\delta_1,$ $\delta_2$ and $\Ggot$ has degree two, we infer that $\Ggot|_{\delta_j}:\delta_j \to \s^2 \cap \{x_3=0\}$ is a diffeomorphism, $j=1,2,$ and  $\Ggot^{-1}(\s^2 \cap \{x_3=0\})=B.$ Furthermore, the planar curve $X(\delta_j)$ is convex, $j=1,2.$

Let $M_a^+$ denote the region $X^{-1}(\{x_3\geq 0\}),$ and note that  $\partial(M_a^+)=B$ and $X|_{M_a^+}$ has a unique catenoidal end. In particular, $\Ggot({M_a^+} )$ is contained in a hemisphere bounded by $\s^2 \cap \{x_3=0\}.$ This fact and the geometry of harmonic catenoidal ends (see Lemma \ref{lem:embe}) imply that  $X(M_a^+)$ is a graph over a region of $\{x_3=0\}$ and  $X:M_a^+ \to X(M_a^+)$ is one to one. To finish, just notice that $X(M)=X(M_a^+) \cup \Tcal(X(M_a^+))$ and use a symmetric argument to infer that $X$ is an embedding.
\end{proof}


\section{The Gaussian image of complete harmonic immersions} \label{sec:gauss}

The questions considered in this section are all related to the size of the Gaussian image of complete harmonic immersions in $\r^3.$

Let us start by proving some Bernstein type results. Theorems \ref{th:mediohuevo} and \ref{th:gaussacotada} below were proved by Klotz \cite{K} under the extra hypothesis that the Hopf differential of the immersion never vanishes. 
\begin{theorem} \label{th:mediohuevo}
Let $M$ be an open Riemann surface, and let $X:M \to \r^3$ be a complete harmonic immersion whose Gaussian image is contained in $\s^2\cap \{x_3>\epsilon\}$ for some $\epsilon>0.$

Then $X(M)$ is a plane.
\end{theorem}
\begin{proof} As usual, label $\Phi$ and $\Ggot$ as the Weierstrass data and the Gauss map of $X,$ respectively. Take a neighborhood $U$  of $(0,0,1)$   in $\s^2$ so that $\Ggot(M) \subset \s^2\cap \{\langle \nu,(x_1,x_2,x_3) \rangle>\epsilon/2\}$ for all $\nu \in U.$ 

For each $\nu \in U,$ denote $\Pi_\nu\subset \r^3$ as the vectorial plane orthogonal to $\nu,$ and write $\pi_\nu:\r^3\to \Pi_\nu$ for the orthogonal projection. Since $\Ggot(M) \subset \s^2\cap \{\langle \nu,(x_1,x_2,x_3) \rangle>\epsilon/2\}$ and $X$ is complete,  the mapping $\pi_\nu \circ X:M \to \Pi_\nu$ is a local diffeomorphism satisfying the path-lifting property. As a consequence,  $\pi_\nu \circ X$ is a global harmonic diffeomorphism, and by Heinz theorem \cite{H} $M$ is biholomorphic to $\c$ (in the sequel, $M=\c$). Let $\{e_1(\nu),e_2(\nu)\}$ be an orthonormal basis of $\Pi_\nu,$ and call $F_j(\nu)=\langle e_j(\nu),\Phi \rangle,$ $j=1,2.$  Since  $\Re\left(F_1(\nu),F_2(\nu)\right):\c \to \r^2$
 is a harmonic diffeomorphism, then $F_j(\nu):\c \to \c$ is a biholomorphism, $j=1,2.$ We infer that  $F_j(\nu)(z)=a_j(\nu) z+b_j(\nu)$ for some constants $a_j(\nu),$ $b_j(\nu) \in \c,$ $a_j(\nu) \neq 0,$   $j=1,2.$  
 
As this holds for any $\nu \in U,$ then $\Phi=\partial_z X=(a_1,a_2,a_3)dz$ for some constants $a_1,$ $a_2$ and $a_3 \in \c,$ so $X(M)$ lies in a plane and we are done. 
\end{proof}

Notice that Theorem \ref{th:mediohuevo} is sharp in the sense that, up to a rigid motion, the Gaussian image of the rotational harmonic horn is $\s^2\cap \{0<x_3<1\}.$

In the quasiconformal case a little more can be said. For instance, it is well known that entire (non necessarily harmonic) QM graphs in $\r^3$ are planes, see \cite{S}. By Lemma \ref{lem:quasiconformal}-(ii),  QC harmonic parameterizations of entire graphs are planes.

In this context, we can also prove the following Osserman type result.

\begin{theorem} \label{th:gaussacotada}
Let $M$ be a non compact Riemann surface with compact boundary, and let $X:M \to \r^3$ be a complete QC harmonic immersion. Assume that the Gaussian image of $X$ lies in $\s^2 \cap \{|x_3|<1-\epsilon\}$ for some $\epsilon >0.$

Then $X$ is of FTC.  Furthermore, if in addition $\partial(M)=\emptyset$  then $X(M)$ is a plane. 
\end{theorem}
\begin{proof}
Let $\Phi=(\Phi_j)_{j=1,2,3}$ and $\Hgot$ denote the Weierstrass data and the Hopf differential of $X.$ Let $g:M \to \overline{\c}$  denote the complex Gauss map of $X,$ and write $2\lambda=\Phi_3 (|g|+1/|g|).$ Notice that from our hypothesis the function $|g|+1/|g|$ is  bounded and $\Phi_3$ never vanishes (see Lemma \ref{lem:gauss}-(a)), hence the metrics $|\lambda|^2$ and $|\Phi_3|^2$ are equivalent. Since $X$ is QC, equation \eqref{eq:beltrami} gives that 
$|1-\Hgot/\lambda^2|<k$ for some constant $k>0.$ Furthermore, by equation \eqref{eq:klotz-g}  ones has that
\[
1\leq \|\Phi\|^2/|\lambda|^2=1+|1-\Hgot/\lambda^2|\leq 1+ k.
\]
Therefore,  the metrics $\|\Phi\|^2,$ $|\lambda|^2$ and $|\Phi_3|^2$ are equivalent. Since $X$ is complete then so is the Klotz metric $\|\Phi\|^2$ (see Remark  \ref{re:kcompleta}). Therefore,   $|\Phi_3|^2$ is a complete flat conformal metric on $M$ and Huber theorem \cite{huber} implies that $M$ is biholomorphic to $\hat{M}-\{P_1,\ldots,P_k\},$ where $\hat{M}$ is a compact Riemann surface with $\partial(\hat{M})=\partial(M)$ and $\{P_1,\ldots,P_k\}\subset \hat{M}-\partial(\hat{M}).$ Furthermore, by Osserman results \cite{O}, the 1-form $\Phi_3$ has a pole at $P_j$ for all $j.$ As $\|\Phi\|^2/|\Phi_3|^2$ is bounded, then $\Phi_j$ is meromorphic on $\hat{M}$ as well, $j=1,2,$ and $X$ is algebraic. By Theorem \ref{th:fun}, $X$ is of FTC. 

Assume in addition that $\partial(M)=\emptyset.$ Since $g$ is quasiconformal, the big Picard theorem for quasiconformal mappings implies that it is constant and $X(M)$ is a plane.
\end{proof}

To finish this section, we state some Privalov type results for harmonic immersions. In other words, we show geometrical conditions for a harmonic immersion to be conformal.

\begin{theorem} \label{th:privalov}
Let $X:\d \to \r^3$ be a non-flat QC harmonic immersion with Weierstrass data $\Phi=(\Phi_j)_{j=1,2,3},$ and label $g$  and $\mu=\partial_{\overline{z}}g/\partial_z g$  as its  complex Gauss map and the Beltrami differential of $g,$ respectively. Let $\Gamma \subset \partial (\d)$ be an open arc, and assume that either of the following conditions holds:
\begin{enumerate}[{\rm (i)}]
\item  The Gaussian image of $X$ lies in $\s^2 \cap \{|x_3|<1-\epsilon\}$ for some $\epsilon>0,$ and $|\mu|$ has angular limit zero almost everywhere on $\Gamma.$ 
\item  $\overline{\c}-g(\d)$ has positive logarithmic capacity, $\lim_{R \to 0} \nu_R=0$ and $\int_{0}^1 \nu_R^2 R^{-2} dR<+\infty,$ where $\nu_R=\sup\{|\mu(x)|\;|\; x \in \Gamma_R\}$ and $\Gamma_R$ is the set of points in $\d$ at distance $R>0$ from $\Gamma.$
\end{enumerate}

Then $X$ is  conformal and incomplete. 
\end{theorem}
\begin{proof} As usual, call $\Hgot$ as the Hopf differential of $X.$

Assume that (i) holds.  Since  $1/|g|+|g|$ is bounded and $|\mu|$ has angular limit zero almost everywhere on $\Gamma,$  then $\Phi_3$ never vanishes and  Proposition \ref{pro:wei-h}-(W.2) implies that the holomorphic function $\Hgot/\Phi_3^2$ has angular limit zero almost everywhere on $\Gamma$ as well. The Luzin-Privalov theorem yields that  $\Hgot=0$ and we are done.

Suppose now that (ii) holds. Set $\hat{\mu}:\overline{\c} \to \c$ by $$\hat{\mu}|_{\d}=\mu,\quad \hat{\mu}|_{\partial(\d)}=0\quad \text{and}\quad \hat{\mu}|_{\overline{\c}-\overline{\d}}=\mu \circ J,$$ where $J:\overline{\c} \to \overline{\c}$ is the antiholomorphic involution $J(z)=1/\overline{z}.$  Let $q:\overline{\c} \to \overline{\c}$ be the unique quasiconformal homeomorphism fixing $0,$ $1$ and $\infty$ and satisfying that $\partial_{\overline{z}} q/\partial_z q= \hat{\mu}.$ We have that $q(\d)=\d$ and $f=g\circ q^{-1}:\d \to \c$ is meromorphic, see \cite{ahlfors}.

\begin{claim} $g$ has angular limits different from both zero and $\infty$ almost everywhere on $\Gamma.$
\end{claim}
\begin{proof} From (ii) one has that $\hat{\mu}$ is continuous on $\Gamma,$ hence $q$ is asymptotically homogeneous at any point $\xi\in \Gamma$ (see \cite{Gutlyanskii-Ryazanov}), that is to say,
\begin{equation} \label{eq:fin}
\lim_{z \to 0} \frac{q(\theta z+\xi)-q(\xi)}{q(z+\xi)-q(\xi)}=\theta\quad \text{for all $\xi \in \Gamma$ and $\theta \in \partial(\d).$}
\end{equation} 
On the hand, since $q$ fixes the unit circle and preserves the orientation one has that  $$\lim_{t \to +0} \frac{q(e^{\imath t}-1+\xi)-q(\xi)}{|q(e^{\imath t}-1+\xi)-q(\xi)|}=\imath q(\xi).$$  Taking into account \eqref{eq:fin},
$$\lim_{t \to +0} \frac{q(\theta (e^{\imath t}-1)+\xi)-q(\xi)}{|q(e^{\imath t}-1+\xi)-q(\xi)|}=\lim_{t\to +0} \frac{q(\theta (e^{\imath t}-1)+\xi)-q(\xi)}{q(e^{\imath t}-1+\xi)-q(\xi)} \;\frac{q(e^{\imath t}-1+\xi)-q(\xi)}{|q(e^{\imath t}-1+\xi)-q(\xi)|}=\imath q(\xi)\theta $$ for any $\theta \in \partial(\d).$

This shows  that $q$  maps non-tangential curves at $\xi$ into non-tangential curves at $q(\xi),$ and likewise for $q^{-1}.$ Furthermore,  the quadratic Carleson measure condition $\int_{0}^1 \nu_R^2 R^{-2} dR<+\infty$ guarantees  that $q$ and $q^{-1}$ have weak derivatives on $\Gamma$ (see \cite{macmanus}), hence they  preserve measure zero sets.
 
On the other hand, as $\overline{\c}-f(\d)$ has positive logarithmic capacity then $f$ has bounded characteristic in the Nevalinna sense (see \cite[p. 213]{nevalinna}), and in particular, $f$ has angular limits almost everywhere on $\partial(\d).$ By the Luzin-Privalov theorem, the angular limits of $f$ are neither $0$ nor $\infty$  on a full subset $E$ of $q(\Gamma)$ (recall that $f$ is not constant). Since $q^{-1}$ preserves measure zero sets, then $q^{-1}(E)$ is a full subset of $\Gamma,$ and so, $g$ has angular limits different from both $0$ and $\infty$ everywhere on $\Gamma.$     
\end{proof}

The above claim, Proposition \ref{pro:wei-h}-(W.2) and the fact $\lim_{R \to 0} \nu_R=0$ imply that $\Hgot/\Phi_3^2$ has angular limit zero almost everywhere on $\Gamma.$ Once again, the Luzin-Privalov theorem yields   $\Hgot=0,$ hence $X$ is conformal. 

Finally, if $X$ were complete  Osserman theorem \cite{O} would imply that $X(\d)$ is a plane,  contradicting Heinz theorem \cite{H}. 
\end{proof}

\begin{remark}
The conclusion of the previous theorem holds with a similar proof  when we replace the hypothesis {\em "$\overline{\c}-g(\d)$ has positive logarithmic capacity"} in {\rm(ii)} for the weaker one {\em "$g$ has angular limits almost everywhere on $\Gamma$". }
\end{remark}

\begin{corollary} \label{co:privalov}
Let $X:M \to \r^3$ be a non-flat AC harmonic immersion, where $M$ is an open Riemann surface carrying  non-constant bounded holomorphic functions. Assume that the  Gaussian image of $X$ lies in $\s^2 \cap \{ |x_3|<1-\epsilon\}$ for some $\epsilon>0.$

Then $X$ is conformal and incomplete.
\end{corollary}
\begin{proof} Fix a non-constant bounded holomorphic function $f:M \to \c.$ Label $\Phi,$ $g$ and $\mu$ as the Weierstrass data  of $X,$ the complex Gauss map of $X$ and the Beltrami differential of $g,$ respectively. 
Let $\pi:\d \to M$ be the conformal universal covering mapping of $M,$ and write $\hat{f}$  and $\hat{\mu}$ for the pull back of $f$  and $\mu$ by $\pi,$ respectively. 

Fatou theorem gives that $\hat{f}$ has (finite) angular limits  on a full set $E \subset \partial(\d).$ 

By Theorem \ref{th:privalov}, we only have to prove that $\hat{\mu}$ has angular limit zero at any $\xi \in E.$ To do this, it suffices to check that for any non tangential Jordan arc $\hat{\alpha}:[0,1) \to \d$ such that $\lim_{t \to 1} \hat{\alpha}(t)=\xi,$  one has  $\lim_{t \to 1} \hat{\mu}(\hat{\alpha}(t))=0.$ Here, ``non tangential'' means that $\hat{\alpha}([0,1))$  lies in a triangle $T \subset \d$ with exterior vertex $\xi.$

Indeed, set $\alpha =\pi\circ \hat{\alpha}$ and call $A\subset M$ as the limit set of $\alpha.$ To be more precise, $P \in A$ if and only if  there exists a  sequence $\{t_n\}_{n \in \n} \subset [0,1)$ with $\{t_n\}\to1$ and $\{\alpha(t_n)\}_{n \in \n} \to P.$  Assume for a moment that $A \neq \emptyset,$ and take a point  $P \in A$ and a closed disc $D\subset M$ with  $P\in D^\circ.$ Label $m$ as the (possibly infinite) number of sheets of $\pi$ and write $\pi^{-1}(D)=\cup_{1\leq j<m+1} D_j,$ where $\pi|_{D_j}:D_j \to D$ is a homeomorphism for all $j$ and $D_{j_1} \cap D_{j_2}=\emptyset$ when $j_1 \neq j_2.$ As $\hat{\alpha}$ is divergent and $P$ is a limit point of $\alpha,$ then necessarily $m=\infty$ and $\pi^{-1}\left(\alpha([0,1))\cap D\right)=\cup_{j=1}^\infty \left(\hat{\alpha}([0,1)) \cap D_j\right)$ contains infinitely many pairwise disjoint  Jordan arcs whose images by $\pi$ accumulate at a continuum $I$ in $A \cap D$  connecting $P$ and $\partial(D).$ However, $f|_I$ is constant and equal to the angular limit of $f$ at $\xi,$ contradicting the identity principle. This shows that $A=\emptyset$ and  $\alpha$ is divergent. The AC hypothesis for  $X$ gives that $0=\lim_{t \to 1} \mu(\alpha(t))=\lim_{t \to 1} \hat{\mu}(\hat{\alpha}(t)).$  

Finally, if $X$ were complete  Osserman theorem \cite{O} would imply that $X$ is flat, a contradiction.
\end{proof}

\end{document}